\documentclass[11pt]{article}
\usepackage[margin=4.46cm]{geometry}
\usepackage{graphicx}
\usepackage{tikz}
\usepgflibrary{shapes}
\usepackage{epsfig}
\usepackage{amsmath}
\usepackage{amssymb}

\usepackage{enumitem}
\usepackage{lipsum}

\usepackage{caption}
\usepackage{colortbl} 
\usepackage{caption}
\usepackage[makeroom]{cancel}
\usepackage{tikz}
\newcommand{\commentout}[1]{}

\def \Rset {{\mathbb R}}

\def \Zset {{\mathbb Z}}

\def \Nset {{\mathbb N}}

\newcommand{\nit}{\noindent}

\newcommand{\be}{\begin{equation}}
\newcommand{\ee}{\end{equation}}
\newcommand{\ba}{\begin{eqnarray}}
\newcommand{\ea}{\end{eqnarray}}
\newcommand{\bi}{\begin{itemize}}
\newcommand{\ei}{\end{itemize}}
\newcommand{\br}{\begin{eqnarray}}
\newcommand{\er}{\end{eqnarray}}

\newcommand{\qed}{\mbox{$\square$}\newline}

\newcommand{\norm}[1]{\left|#1\right|}

\newcommand{\Div}[1]{\text{div}\left(#1\right)}

\newcommand{\inner}[2]{\left\langle#1,#2\right\rangle}
\newcommand{\Dpn}[1]{\frac{D^\perp #1}{\left|D #1\right|}}

\newtheorem{theo}{Theorem}[section]
\newtheorem{defin}{Definition}[section]

\newtheorem{lem}{Lemma}[section]
\newtheorem{cor}{Corollary}[section]
\newtheorem{rmk}{Remark}[section]

\numberwithin{equation}{section}

\begin{document}
\title{Existence of {an effective burning velocity in a cellular flow for the curvature G-equation 
proved using a game analysis}}
\author{Hongwei Gao \thanks{Morning Star Institute, Irvine, CA, Email: hwgao2@gmail.com}, \ Ziang Long 
\thanks{Meta Inc, Mountain View, CA. Email: zlong6@uci.edu}, \ Jack Xin, \ Yifeng Yu \thanks{Department of Mathematics, University of California at Irvine, Irvine, CA 92697.
Email: (jack.xin,yifengy)@uci.edu. The work was partly supported by 
NSF grants DMS-1952644/DMS-2309520 (JX) and DMS-2000191 (YY).}}

\date{}
\maketitle

\begin{abstract}
G-equation is a popular level set model in turbulent combustion, and 
becomes an advective mean curvature type evolution equation when curvature of a moving flame
in a fluid flow is considered: 
$$
G_t + \left(1-d\, \Div{\frac{DG}{|DG|}}\right)_+|DG|+V(x)\cdot DG=0.
$$
Here $d>0$ is the Markstein number and the positive part $()_+$ is imposed to 
avoid a non-physical negative laminar flame speed.
 {For simplicity of presentation, we focus
mainly on the case when $V:\Rset^2\to \Rset^2$ is the two dimensional cellular 
flow with Hamiltonian $H = \sin x_1 \, \sin x_2$ and amplitude $A$.
Our main result is that} for any unit vector $p\in \Rset^2$, 
there exists a positive number $\overline H(p)$ such that if $G(x,0)=p\cdot x$, then
$$
\left|G(x,t)-p\cdot x+\overline H(p)t\right|\leq C \quad \text{in $\Rset^2\times [0,\infty)$}
$$
for a constant $C$ depending only on 
{on the Markstein number $d$ and the cellular 
flow amplitude $A$.} The number $\overline H(p)$ corresponds to 
the effective burning velocity in the physics literature. 
The non-coercivity encountered here is one of the major difficulties for 
homogenization of the mean curvature-type equations.  To overcome it,  we introduce 
a new approach that combines PDE methods with a dynamical analysis of the 
Kohn-Serfaty deterministic game characterization of the curvature G-equation 
utilizing the streamline structure of cellular flows. Extension to general 
two-dimensional incompressible flows is also discussed.  
In three dimensional incompressible flows, the existence of $\overline H(p)$ might fail when the 
flow intensity exceeds a bifurcation value even for simple shear flows \cite{MMTXY2023}.

\end{abstract}
\vspace{.2 in}

\hspace{.12 in} {\bf AMS Subject Classification:} 70H20, 76M50, 76M45, 76N20.
\bigskip

\hspace{.12 in} {\bf Key Words:} Front motion, cellular flow, curvature effect,

\hspace{.12 in} homogenization, game theory.

\thispagestyle{empty}
\newpage

\section{Introduction}
\setcounter{equation}{0}
\setcounter{page}{1}

G-equation is a well-known level-set model in turbulent combustion \cite{W1985}. 
Precisely speaking, given a reference function $G(x,t)$, let the flame front be the zero level set $\{G(x,t)=0\}$ at time $t$, where 
the burnt and unburnt regions are $\{G(x,t) < 0\}$ and $\{G(x,t) > 0\}$, respectively. 
See Figure \ref{levset} below. The velocity of ambient fluid $V:\Rset^n\to \Rset^n$ is 
assumed to be continuous and $\Zset^n$-periodic. 
The propagation of flame front obeys a simple motion law: 
$$
{v}_{\vec{n}}=s_l+V(x)\cdot \vec{n},
$$ 
i.e., the normal velocity is the laminar flame speed ($s_l$) plus the projection of $V$ 
along the normal direction $\vec{n}$. This leads to the so--called $G$-equation, a level-set 
PDE \cite{W1985,P2000,OF2002}:
$$
G_t + V(x)\cdot DG + s_l \, |DG|=0 \quad \text{in $\Rset^n\times (0,+\infty)$}. 
$$

\begin{center}\label{levset}
\includegraphics[scale=0.8]{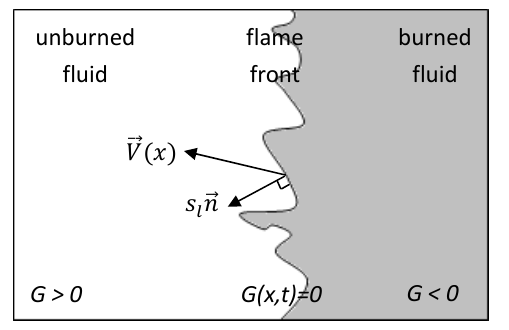} \quad \qquad
\captionof{figure}{Level-set formulation of front propagation.}
\end{center}

\bigskip

In the simplest case when the laminar flame speed $s_l$ is a positive constant, the above 
equation is called the inviscid G-equation, a non-coercive convex Hamilton-Jacobi equation. 

However, in general, $s_l$ is not constant along the flame front since the 
burning temperature varies in different locations. In order to quantify 
temperature difference, the curvature effect in turbulent combustion was first 
introduced by Markstein \cite{Mar}, which says that if the flame front bends 
toward the cold region (unburned area, point C in Figure 2 below), 
the flame propagation slows down. If the flame front bends toward the hot 
spot (burned area, point B in Figure 2), it burns faster.
\begin{center}\label{curv}
\includegraphics[scale=0.6]{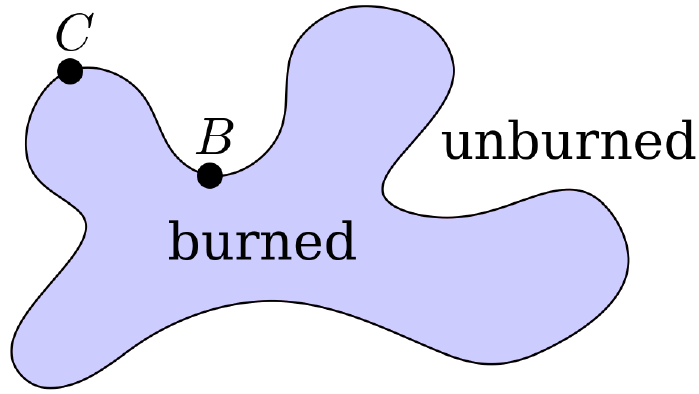}
\captionof{figure}{Curvature effect.}
\end{center}

There is a {vast turbulent combustion literature} discussing the impact of curvature effect. 
Below is the most recognized empirical linear relation proposed by Markstein \cite{Mar} to 
approximate the dependence of the laminar flame speed $s_l$ on the curvature (see also \cite{P2000}):
\be\label{laminar}
s_l=s_{l}^{0}(1- d\; \kappa)_{+}. 
\ee
Here $s_{l}^{0}$, the mean value, is a positive constant. Hereafter we set $s_{l}^{0}=1$. The parameter 
$d>0$ is the so called Markstein length which is proportional to the flame thickness. 
The mean curvature along the flame front is $\kappa$. The positive 
part $(a)_{+}=\max\{a,0\}$ is imposed to avoid negative laminar flame speed since materials 
cannot be ``unburned". This correction usually is not explicitly mentioned in combustion literature since, 
by default, the curvature is always assumed to be small there. However, mathematically, large positive curvature could occur as time evolves. Therefore it is necessary to 
explicitly add this correction in theoretical study and numerical computations \cite{ZR_94} 
if the physical validity is taken into consideration in the modeling of flame propagation. 

Plugging the expression of the laminar flame speed (\ref{laminar}) 
into the G-equation and normalizing the constant $s_{l}^{0}=1$, we obtain a mean curvature type 
equation with advection:
\be
G_t + \left(1- d\, \, \mathrm{div}\left({DG\over |DG|}\right)\right)_{+}|DG|+V(x)\cdot DG=0.\label{ge1}
\ee
We would like to mention that curvature G-equations served as some of the first 
numerical examples in the development of level-set methods \cite{Se1985,OS1988,OF2002}. 
The rigorous mathematical foundation under the framework of viscosity solution was 
established in \cite{CGG, ES}.

In general, $\kappa$ changes sign along a curved flame front. 
So a mathematically interesting and physically important question is:

\medskip 

\nit {\bf Question 1:} {\it How does the ``averaged" flame propagation speed depend on the curvature term?}

\medskip

As the first step, the main purpose of this paper is to establish the existence of a properly 
defined ``averaged speed", which is basically to average fluctuations caused by both 
the flow and the curvature. The theory of homogenization provides such a rigorous 
mathematical framework in environments with microstructures.

Turbulent combustion usually involves small scales. As a simplified model, 
we rescale $V$ as $V=V({x\over \epsilon})$ and replace $d$ by $d\, \epsilon$. 
Here $\epsilon$ denotes the Kolmogorov scale (the small scale in the flow). 
The diffusivity constant $d>0$ is called the Markstein number. 
We would like to point out that the dimensionless Markstein number is 
$d\cdot {\delta_L \over \epsilon}$ with $\delta_L$ denoting the flame thickness \cite{P2000}. 
In the thin reaction zone regime, $\delta_L = O(\epsilon)$, see Eq. (2.28) and Fig. 2.8 of \cite{P2000}. 
Without loss of generality, let ${\delta_L \over \epsilon}=1$. Then (\ref{ge1}) becomes
\be
{\partial G_{\epsilon}\over \partial t }+\left(1-d\, \epsilon \, \mathrm{div}\left({DG_{\epsilon}\over |DG_{\epsilon}|}\right)\right)_{+} |DG_{\epsilon}|+ V\left({x\over \epsilon}\right)\cdot DG_{\epsilon}=0.\label{ge1b}
\ee

Since $\epsilon \ll 1$, it is natural to look at $\lim_{\epsilon\to 0}G_{\epsilon}$, i.e., 
the homogenization limit. If $G_{\epsilon}(x,0)=p\cdot x$, the limit $\lim_{\epsilon\to 0}G_{\epsilon}(x,t)$ can be viewed as the flame propagation under the effective burning velocity (also called {\it ``turbulent flame speed"} in physics literature) along a given direction $p$.

For general $V$, the effective burning velocity might not exist (see Remark \ref{rmk:nonexist}). 
Moreover, a general $V$ is unlikely to provide meaningful answers to 
delicate problems like Question 1.  
In this paper, we consider 
a two-dimensional (2D) cellular flow that is frequently used in the mathematics 
and physics literature. See \cite{CG} for more physical motivations. 
For clarity of presentation, we choose to work with the typical example
$$
V(x)=A(DH)^{\perp}=A(-\cos x_2\sin x_1,\ \cos x_1\sin x_2).
$$
Here $H(x_1,x_2)=\sin x_1\sin x_2$ is the stream function and the positive constant $A$ represents the flow intensity. The following is our main theorem.

\begin{theo}\label{theo:main} Suppose that $V=A(DH)^{\perp}$ for $A>0$. For each $\epsilon>0$, let 
$G_{\epsilon}=G_{\epsilon}(x,t)\in C(\Rset^2\times [0, \infty))$ be the unique viscosity solution to equation (\ref{ge1b}) on $\Rset^2\times [0, \infty)$ subject to $G_{\epsilon}(x,0)=p\cdot x$ for a given $p\in \Rset^2$. Then
$$
|G_{\epsilon}(x,t)-p\cdot x+\overline H_A(p)t|\leq C\epsilon \quad \text{on $\Rset^2\times [0,\infty)$}
$$
for a constant $C$ depending only on $d$, $V$ and $|p|$. Here $\overline H_A\in C(\Rset^2\backslash\{0\}, (0,\infty))$ is a continuous positive homogeneous function of degree one. 

Equivalently, if $G(x,t)$ is the unique viscosity solution to equation (\ref{ge1}) on $\Rset^2\times [0, \infty)$ subject to $G(x,0)=p\cdot x$, then
$$
\left |G(x,t)-p\cdot x+\overline H_A(p)t\right|\leq C \quad \text{on $\Rset^2\times [0,\infty)$}.
$$

\end{theo}

\begin{rmk}\label{rmk:fullhom} Our proof relies on establishing the existence of 
correctors to the cell problem for every $p\in \Rset^2$:
$$
\left(1- d\, \, \mathrm{div}\left({p+Dv\over |p+Dv|}\right)\right)_{+}|p+Dv|+V(x)\cdot (p+Dv)=\overline H(p).
$$
In general, the cell problem is sufficient to obtain the full homogenization with any uniformly continuous initial data by Evans' 
perturbed test function method \cite{Evans}. However, as it was first explicitly 
pointed out in \cite{CM}, due to the discontinuity of the mean curvature operator, 
the standard perturbed test function method is not directly applicable. Roughly speaking, given a test function $\phi(x,t)\in C^{\infty}(\Rset^n\times(0,\infty))$, let $p=D_x\phi(x_0,t_0)$, it is not clear whether 
$$
\phi(x,t)+\epsilon v\left({x\over \epsilon}\right)
$$
is an approximate solution of equation (\ref{ge1b}) near $(x_0,t_0)$ due to possible vanishing gradients. 
The authors of \cite{CM} have introduced a mechanism to implement a modified version of 
the perturbed test function method. { See also \cite{AC2018} for a different modified argument  when the corresponding cell problem has a classical solution. To avoid a lengthy paper, we will
address this subtle issue in a future publication.
In view of our problem's physical origin, it would be
an interesting project } to study more detailed properties of 
$\overline H_A(p)$ (e.g., its anisotropy). See Remark \ref{rmk:property}. 
\end{rmk}

\begin{rmk}We would like to point out that the small scales in the curvature G-equation 
help flame propagation. In fact, 
it is not hard to see that under equation (\ref{ge1}), {\it the flame front is stagnant if 
the initial burned region is $\{x\in \Rset^2|\ G(x,0)<0\}=\{x\in [0, \pi]\times [0, \pi]|\ H(x)>1-r\}$ and the initial flame front is 
$$
\{x\in \Rset^2|\ G(x,0)=0\}=\{x\in [0, \pi]\times [0, \pi]|\ H(x)=1-r\}
$$
when $r$ is small enough}. Nevertheless, {\it after homogenization, 
the front moves forward with a positive normal speed $\overline H_A(p)$ 
along any given direction $p$}. The curvature term plays a subtle role of flame spreading when 
the cellular flow corrugates the level set. 
Such a flow induced geometric mechanism is absent in the invicid equation model (i.e. $d=0$). 
The homogenization of the inviscid G-equation has been established 
independently in \cite{CNS} and \cite{XY2010} via completely different approaches 
for general incompressible flows in any dimension.

\end{rmk}

\begin{rmk}
{As pointed out by a referee,  a natural question is whether the cut-off $(\cdot)_{+}$ in the governing equation 
is actually active (equal to zero) at some point $(x,t)$ 
during the evolution in the cellular flow.  A direct numerical simulation suggests that 
the answer is positive 
over $A \in (0, 45]$ at $d=0.1$ where the $\overline{H}$ differs with and without the cutoff. That the cut-off is active is clear in the case of three dimensional shear flows
 \cite{MMTXY2023} where the homogenization fails for large flow intensity 
with  cut-off $(\cdot)_{+}$ and the homogenization holds for 
all flow intensity without cut-off $(\cdot)_{+}$. See Remark \ref{rmk:nonexist}  for more details. }
\end{rmk}

The major remaining question is to understand how the turbulent velocity $\overline H_A(p)$ depends on the Markstein number $d$ and the flow intensity $A$. There is a consensus in combustion literature that the curvature effect slows down flame propagation \cite{Ron}. 
Heuristically, this is because the curvature term smoothes out the flame front and reduces the total area of chemical reaction \cite{Se1985}. To the best of our knowledge, the first mathematically rigorous result in this direction was obtained in \cite{LXY} for shear flows. It was proved that $\overline H$ is strictly decreasing with respect to the Markstein number $d$, which is consistent with the experimental observation (e.g., \cite{CWL2013}).

For the more complicated cellular flow, it is probably more convenient to compare the growth law with respect to the flow intensity $A$. For the inviscid case, a sharp growth law was establish in \cite{XY2013}: when $d=0$, for $p=(p_1, p_2)$, 
$$
{A\, \pi(|p_1|+|p_2|)\over 2\log A+C_2}\leq \overline H_A(p) \leq {A\, \pi(|p_1|+|p_2|)\over 2\log A+C_1}.
$$

By refining the method and estimates  in this paper, we conjecture  that 
{it might be possible to show} that when $d>0$, 
$$
\overline H_A(p)=O\left({A\over \log A}\right).
$$
{More interestingly, in future work we plan to investigate:}
\medskip

{\bf Question 2:} Does the presence of the curvature effect significantly reduce the prediction of the turbulent flame speed? More specifically, does there exist a unit direction $p$ such that 
$$
\limsup_{A\to +\infty}{\overline H_A(p)\log A\over A}<{\pi(|p_1|+|p_2|)\over 2}\ ?
$$

If the curvature term is replaced by full diffusion (i.e. the Laplacian $\Delta$ without any correction), a 
dramatic slow-down ($\overline H_A(p)$ uniformly bounded in $A$) is proved in \cite{LXY2011} 
for two dimensional (2D) cellular flows.

\medskip

{\bf $\bullet$ Other works.} In non-combustion practical suitations such as phase transition (e.g. crystal growth) in material science \cite{L1980,CB},  the  motion law is given by (no drift term and $()_+$ correction)
\be\label{eq:GT-relation}
v_{\vec{n}}=a(x)-d\kappa
\ee
for a continuous positive $\Zset^n$-periodic function $a(x)$. The above formula is known as 
the Gibbs-Thomson relation. The $()_+$ is not needed in crystal growth since both 
freezing and melting could occur in the 
situation of ice formation \cite{L1980}.   See (\cite{LS2005, DKY, CLS, CB, CM, CN, KG, GMT,AC2018, Feldman}, etc) and reference therein for works regarding the homogenization or large time limit related to (\ref{eq:GT-relation}).   When $a(x)$ satisfies a special coercivity condition, 
homogenization has been proved in \cite{LS2005} for all dimensions. 
When $n=2$, homogenization was proved in \cite{CM} for all positive $a(x)$ by  a geometric approach.  Moreover, 
the authors have also constructed a counterexample in \cite{CM} 
when the dimension $n\geq 3$ for a positive $a(x)$. 
See \cite{C2013} for more discussions about geometric approaches in homogenization.  
 Due to the presence of the non-coercive transport term $V\cdot DG$, 
{or $a= a(x, \vec{n})= 1 + V(x) \cdot \vec{n}$},  
our situation and methods are very different from all these previous works. See Remark \ref{rmk:cm} for more discussions.  
Our result also holds when $(1-d\kappa)_+$ is replaced by a more general form $(a(x)-d\kappa)_+$.  

\medskip

\begin{rmk}\label{rmk:nonexist} Like the invisicd G-equation case ($d=0$), 
it is easy to construct a smooth periodic $V$ such that the effective speed 
does not exist, e.g. let $V(x)=K(x-P_0)$  near $P_0=({\pi\over 2}, {\pi\over 2})$ for a 
large number $K>0$ and equal to the cellular flow away from $P_0$.  
We conjecture that our method could be  extended to cover  general 
2D smooth periodic incompressible flows ($\mathrm{div} (V)=0$). 
See Section \ref{sec:extension} for more detailed discussions.  
When $n\geq 3$,  it was proved in \cite{MMTXY2023} that, 
for the shear flow $V(x_1,x_2,x_3)=(0,0,Af(x_1,x_2))$, 
the effective burning velocity ceases to exist when the flow 
intensity $A$ exceeds a bifurcation value. Note that when $n=2$, 
the existence of effective burning velocity can be established 
very easily for all shear flows since the cell problem is reduced 
to an ODE and maximum principle immediately leads 
to uniform bound of derivatives. 
\medskip

{\bf Question 3:} Does Theorem \ref{theo:main} hold if $V$ is a 
3D physically meaningful flow that possesses turbulent or {chaotic}  
structures, e.g. the Arnold-Beltrami-Childress flow \cite{CG,XYZ}?
\end{rmk}

\medskip

{\bf $\bullet$ Outline of the paper.} In section 2, for reader's convenience, 
we will go over the definition and comparison principle of viscosity solutions of 
curvature type equations. 
The associated game theoretic interpretation from \cite{KS1} will be reviewed as well. 
In section 3, we prove Theorem \ref{theo:main} by a novel approach combining 
Lagrangian (game dynamics) and Eulerian approaches (PDE techniques).  
In section 4, we discuss how to extend our results 
and methods to general 2D incompressible flows. 
Our approach also suggests a possible {\it general framework} to tackle  
non-periodic settings (see Remark \ref{rmk:non-periodic}).  
{To help derive a reachability property in our proof, 
we show rigorously in the Appendix  
a useful fact related to the consistency of viscosity and classical solutions in front propagation.}

\section{Preliminary}

In this section, $V:\Rset^n\to \Rset^n$ is assumed to be continuous and $\Zset^n$-periodic. Write
$$
F(A,p)=\left(|p|-d\left(\mathrm{tr}A-{p\cdot A\cdot p\over |p|^2}\right)\right)_{+}
$$
for $(A,p)\in S^{n\times n}\times \left(\Rset^n\backslash\{0\}\right)$. Here $S^{n\times n}$ is the space of $n\times n$ real symmetric matrices. Following the definition in \cite{CL1992}, let 
$$
\underline{F}(A,p)=
\begin{cases}
F(A,p) \quad \text{if $|p|\not=0$}\\
0\quad \text{if $|p|=0$}
\end{cases}
\quad \mathrm{and} \quad 
\overline{F}(A,p)=
\begin{cases}
F(A,p) \quad \text{if $|p|\not=0$}\\
2dn||A|| \quad \text{if $|p|=0$}.
\end{cases}
$$
Here $||A||$ is the largest absolute value of eigenvalues of $A$. Note that $\underline F$ and $\overline F$ are lower semicontinuous and upper semicontinuous respectively.  As it was pointed out in \cite{CL1992}, the definition at $p=0$ is not really important as long as the lower and upper continuity hold. This is because  $A$ will be  zero when $p=0$ with the proper test function (see the proof of Theorem \ref{theo:comparison} for instance).  When $n=2$,
$$
\mathrm{tr}A-{p\cdot A\cdot p\over |p|^2}={p^{\perp}\cdot A\cdot p^{\perp}\over |p|^2}.
$$
Here $p^{\perp}=(-p_2,p_1)$ if $p=(p_1,p_2)\in \Rset^2\backslash\{0\}$. 

\subsection{Definition of viscosity solutions and comparison principle}

We first introduce several definitions and terminologies. 

Let
$$
\text{$LSC(\Omega)$: the set of lower semicontinuous functions defined on $\Omega$},
$$
$$
\text{$USC(\Omega)$: the set of upper semicontinuous functions defined on $\Omega$},
$$
and
$$
\text{$UC(\Omega)$: the set of uniformly continuous functions defined on $\Omega$}.
$$

\begin{defin} Assume that $G=G(x,t)\in USC(\Rset^n\times (0, \infty))$. $G$ is called a viscosity subsolution of equation (\ref{ge1}) provided that given $\phi(x,t)\in C^2(\Rset^n\times (0, \infty))$, if 
$$
G(x_0,t_0)-\phi(x_0,t_0)=\max_{\Rset^n\times (0,\infty)} (G(x,t)-\phi(x,t)),
$$
then
$$
\phi_t(x_0,t_0)+\underline{F}\left(D^2\phi(x_0,t_0), D\phi (x_0,t_0) \right)+V(x_0)\cdot D\phi(x_0,t_0)\leq 0.
$$
\end{defin}

\begin{defin} Assume that $G=G(x,t)\in LSC(\Rset^n\times (0, \infty))$. $G$ is called a viscosity supersolution of equation (\ref{ge1}) provided that given $\phi(x,t)\in C^2(\Rset^n\times (0, \infty))$, if 
$$
G(x_0,t_0)-\phi(x_0,t_0)=\min_{\Rset^n\times (0,\infty)} (G(x,t)-\phi(x,t)),
$$
then
$$
\phi_t(x_0,t_0)+\overline{F}\left( D^2\phi(x_0,t_0), D\phi (x_0,t_0)\right)+V(x_0)\cdot D\phi(x_0,t_0)\geq 0.
$$
\end{defin}

\begin{defin}$G=G(x,t)\in C(\Rset^n\times (0,\infty))$ is called a viscosity solution of equation (\ref{ge1}) if it is both a viscosity subsolution and a viscosity supersolution. If the initial data $G(x,0)=g(x)$ is given, then we require that
$$
\lim_{t\to 0}G(x,t)=g(x) \quad \text{locally uniformly for $x\in \Rset^n$}. 
$$
\end{defin}

The following comparison principle for solutions can be proved by standard approaches (\cite{CGG, ES, CL1992, CM}), which is a special case of Theorem 3.3 in \cite{CM} . For the reader's convenience, we provide a sketch of the proof. See the proof of Theorem 3.3 in \cite{CM} (the long version)  
for more general operators.

\begin{theo}\label{theo:comparison} Assume that $V\in W^{1,\infty}(\Rset^n)$, $g_1, g_2\in UC(\Rset^n)$ and $g_1\leq g_2$. Suppose that  for some $T>0$, $G_1=G_1(x,t)\in USC(\Rset^n\times (0, T))$ is a viscosity subsolution of equation (\ref{ge1}) and $G_2=G_2(x,t)\in LSC(\Rset^n\times (0, T))$ is a viscosity supersolution of equation (\ref{ge1}). 
If both $G_1$ and $G_2$ are subject to (for $i=1,2$): 
$$
\lim_{t\to 0}|G_i(x,t)-g_i(x)|=0 \quad \text{locally uniformly in $\Rset^n$},
$$
then
$$
G_1(x,t)\leq G_2(x,t) \quad \text{for all $(x,t)\in \Rset^n\times [0, T)$}.
$$

\end{theo}

Proof: By considering ${2\over \pi}\arctan({G_1})$ and ${2\over \pi}\arctan({G_2})$, we may assume that for $i=1,2$, 
$$
-1\leq G_i\leq 1.
$$

We argue by contradiction. If not, then 
$$
\sup_{(x,t)\in \Rset^n\times [0,T)}(G_1(x,t)-G_2(x,t))>0.
$$
Let $M=4(||V||_{W^{1,\infty}(\Rset^n)}+1)$ and 
$$
\tilde G_i=e^{-Mt}G_i.
$$
Then $\tilde G_1$ is a viscosity subsolution of 
$$
{\partial \tilde G_{1}\over \partial t}+M\tilde G_1+ \underline F(D^2\tilde G_1, D\tilde G_1)+V(x)\cdot D\tilde G_1=0
$$
and $\tilde G_2$ is a viscosity supersolution of 
$$
{\partial \tilde G_{2}\over \partial t}+M\tilde G_2+ \overline F(D^2\tilde G_2, D\tilde G_2)+V(x)\cdot D\tilde G_2=0
$$
subject to $ \lim_{t\to 0}|\tilde G_i(x,t)-g_i(x)|=0$ locally uniformly in $\Rset^n$. 

We also have that
$$
r=\sup_{(x,t)\in \Rset^n\times [0,T)}(\tilde G_1(x,t)-\tilde G_2(x,t))>0.
$$
Now let 
$$
w(x,y,t)=\tilde G_1(x,t)-\tilde G_2(y,t)-\Phi_{\rho,K,\delta}(x,y)
$$
for some $K, \rho, \delta>0$ and 
$$
\Phi_{\rho,K, \delta}(x,y,t)=\rho\left(\sqrt{|x|^2+1}+\sqrt{|y|^2+1}\right)+K|x-y|^4+{\delta\over {T-t}}.
$$

{\bf Step 1:} It is easy to see that there exist $\rho_0,\delta_0\in (0,1)$ such that for all $\rho\leq \rho_0$, $\delta\leq \delta_0$ and $K\in \Nset$,
$$
\max_{\Rset^n\times \Rset^n\times [0,T)} w(x,y,t)\geq {1\over 2}\sup_{(x,t)\in \Rset^n\times [0,T)}(\tilde G_1-\tilde G_2)={r\over 2}.
$$
Choose $(\bar x, \bar y, \bar t)\in \Rset^n\times \Rset^n\times [0,T)$ such that
$$
w(\bar x,\bar y,\bar t)=\max_{\Rset^n\times \Rset^n\times [0,T)} w(x,y,t).
$$
For simplicity, we omit the dependence of $(\bar x, \bar y, \bar t)$ on $\rho$, $K$ and $\delta$. 

\medskip

{\bf Step 2 (avoiding the boundary $t=0$)}: Due to $|\tilde G_i|\leq 1$, we have that 
$$
\Phi_{\rho, K,\delta}(\bar x,\bar y,\bar t)\leq 2.
$$
Since for $i=1,2,$ $g_i\in UC(\Rset^n)$, there exists $K_0\in \Nset$ independent of $\rho$ and $\delta$ such that $\bar t>0$ when $K\geq K_0$.

\medskip

{\bf Step 3 (plugging into the equation)}: Now we fix $K=K_0$ and $\delta=\delta_0$. Write $\gamma(x,y)=\rho\left(\sqrt{|x|^2+1}+\sqrt{|y|^2+1}\right)$. Owing to Theorem 8.3 and Remark 3.8 in \cite{CL1992}, there exist $a,b\in \Rset$,  two $n\times n$ symmetric matrices $X$ and $Y$ such that
$$
a-b={\delta\over (T-t)^2}, \quad X\leq Y, \quad ||X||+||Y||\leq CK_0|\bar x-\bar y|^2
$$
and at point $(\bar x, \bar y, \bar t)$
$$
\begin{array}{ll}
a+M\tilde G_1(\bar x,\bar t)+\underline F(X+D_{x}^{2}\gamma, \bar p+D_x\gamma)+V(\bar x)\cdot (\bar p+D_x\gamma)\leq 0\\[3mm]
b+M\tilde G_2(\bar y, \bar t)+\overline F(Y-D_{y}^{2}\gamma, \bar p-D_y\gamma)+V(\bar y)\cdot (\bar p-D_y\gamma)\geq 0
\end{array}
$$
for $\bar p=4K_0(\bar x-\bar y)|\bar x-\bar y|^2$. Since   $|D_{x}^2\gamma|+|D_{y}^2\gamma|+|D_{x}\gamma|+|D_{y}\gamma|\leq C\rho$,
$$
\tilde G_1(\bar x,\bar t)-\tilde G_2(\bar y, \bar t)\geq {r\over 2}+K_0|\bar x-\bar y|^4,
$$
$$
\quad |\bar p|\leq 4K_0|\bar x-\bar y|^3 \quad \mathrm{and} \quad \ |V(\bar x)-V(\bar y)|\leq \left({M\over 4}-1\right)|\bar x-\bar y|,
$$
we derive that
\be\label{F-comparison}
\underline F(X+D_{x}^{2}\gamma, \bar p+D_x\gamma)-\overline F(Y-D_{y}^{2}\gamma, \bar p-D_y\gamma)\leq -{Mr\over 2}-{\delta_0\over (T-\bar t)^2}+C\rho. 
\ee

{\bf Step 4}: Let us look at the following two cases (up to a subsequence if necessary). 

\medskip

Case 1: $\lim_{\rho\to 0}|\bar x-\bar y|=0$. Then both the gradient and the Hessian of the test function go to zero, i.e., 
$$
\begin{array}{ll}
\lim_{\rho\to 0 } \bar p=0\\[3mm]
\lim_{\rho\to  0} X=\lim_{\rho\to 0 } Y=0.
\end{array}
$$
By (\ref{F-comparison}), this implies $0<-{Mr\over 2}$ after sending $\rho\to 0$, which is absurd. 

\medskip

Case 2: $\lim_{\rho\to 0}|\bar x-\bar y|\not= 0$. Without loss of generality, we may assume that 
$$
\lim_{\rho\to 0}\bar p=p_0\not=0, \quad \lim_{\rho\to 0 }X=X_0 \quad \mathrm{and} \quad \lim_{\rho\to 0 } Y=Y_0.
$$
Since $F(A, p)$ is monotonically decreasing with respect to the symmetric matrix $A$ and $X_0\leq Y_0$, taking limit $\rho\to 0$ on (\ref{F-comparison}) leads to 
$$
0\leq F(X_0,p_0)-F(Y_0,p_0)\leq -{Mr\over 2}.
$$
This is impossible. \qed

\begin{cor}\label{cor:lineargrowth} Assume that $g\in UC(\Rset^n)$ and $|g(x)|\leq C(|x|+1)$ for some positive constant $C$. Suppose that $G=G(x,t)\in C(\Rset^n\times [0, \infty))$ is a viscosity subsolution of equation (\ref{ge1}) subject to 
$G(x,0)=g(x)$. Then there exists a constant $\tilde C$ depending only on $C$, $V$ and $d$ such that
$$
|G(x,t)|\leq \tilde C(1+|x|+t) \quad \text{for all $(x,t)\in \Rset^n\times [0, \infty)$}.
$$
\end{cor}

Proof: Let $\phi(x)=2C\sqrt{|x|^2+1}$. Let
$$
C_1=\max_{x\in \Rset^n}\left|\overline F(D^2\phi(x), D\phi(x))+V(x)\cdot D\phi(x)\right|.
$$
Then $G_1(x,t)=\phi(x)+C_1t$ is a viscosity supersolution to equation (\ref{ge1}). Hence by the comparison principle Theorem \ref{theo:comparison}, we have that
$$
G(x,t)\leq G_1(x,t)\leq \tilde C(1+|x|+t)
$$
for $\tilde C=\max (2C, C_1)$. A similar inequality holds in the other direction. 
\qed

\subsection{Game theory interpretation from \cite{KS1}}

In this section, let $n=2$ and $x\in \Rset^2$. Fix the game step size parameter $\tau>0$ and the number of steps $N\in \Nset$. 

For $n\geq $ 1, consider the discrete dynamical system $\{x_n\}_{n=1}^{N}\subset\Rset^2$ associated with the game starting from $x_0=x$: for $n=0,1,2,.., N-1$
$$
\begin{cases}
x_{n+1}=x_n+\tau\sqrt{2d}b_n\eta_n+\tau^2\eta_{n}^{\perp}-\tau^2V(x_n)\\
x_0=x, 
\end{cases}
$$
where $\eta_n\in \overline B_1(0)=\{v\in \Rset^2|\ |v|\leq 1\}$ and $b_n\in \{-1,1\}$. Here for convenience, compared to \cite{KS1}, we reverse the time that leads to $-V$ instead of $V$. Moreover,  $\eta^{\perp}=(-c, a)$ if $\eta=(a,c)\in \Rset^2$.

There are two players: player I and player II. In each step, 

\medskip

$\bullet$ Player I: First choose the direction $\eta_n$;

\medskip

$\bullet$ Player II: Then choose the sign $b_n$. 

\medskip

The sequence of positions $\{x_k\}_{k=0}^{N}$ is called a game trajectory associated with $(\tau, N)$.

\medskip

The goal of player I is to minimize $g(x_N)$ and player II aims to maximize $g(x_N)$ . 

A strategy  $\Gamma_1$ of player I refers to how player I chooses $\eta_n$ based on $\{x,\eta_i, b_i|\ 1\leq i\leq n-1\}$ for $n\geq 1$. Similarly, a strategy $\Gamma_2$ of player II refers to 
how player II chooses the sign $b_n$ based on $\{x, \eta_i, \eta_n, b_i, |\ 1\leq i\leq n-1\}$ for $n\geq 1$.

Assuming that both players play optimally. Write the value function of player I as
\begin{equation}\label{gamedef}
\text{$u_\tau(x,N\tau^2)=$ infimum of $g(x_N)$}.
\end{equation}
Then we have the following dynamic principle: for $k\geq 1$, 
\begin{equation}\label{gamedp}
u_\tau(x,k\tau^2)=\inf_{\norm{\eta}\leq 1}\max_{b=\pm1}u_\tau\left(x+\tau b\sqrt{2d}\eta+\tau^2\eta^\perp-\tau^2V(x),\ (k-1)\tau^2\right).
\end{equation}

\medskip

For the standard mean curvature motion, player I can only choose unit vectors. Here we basically allow player I to also choose $\eta_n=0$, which excludes the possibility of ``unburning". The consistency holds: for small $\tau$,
$$
\begin{array}{ll}
&\min_{\eta\in \overline B_1(0)}\max_{b=\pm 1}\{ \tau b\sqrt{2d}\eta\cdot D\phi (x)+d \tau^2 \eta\cdot D^2\phi(x)\cdot \eta+\tau^2\eta^{\perp}\cdot D\phi\}\\[5mm]
&\in -\tau^2[\underline F(D^2\phi(x), D\phi(x)),\ \overline F(D^2\phi(x), D\phi(x))].
\end{array}
$$

\subsection{Convergence}
Given $g(x)\in UC(\Rset^2)$, as in \cite{KS1}, we have that 
\be\label{limit}
\lim_{\tau\to 0,\ N\tau^2\to t}u_{\tau}(x,N\tau^2)=u(x,t) \quad \text{locally uniformly on $\Rset^2\times (0, \infty)$}.
\ee
Here $u\in C(\Rset^2\times [0,\infty))$ is the unique viscosity solution to equation (\ref{ge1}) subject to 
$$
\lim_{t\to 0}u(x,t)=g(x) \quad \text{uniformly in $\Rset^2$}.
$$
As it was mentioned in \cite{KS2}, this game theory  approach is  essentially a semi-discrete approximation scheme (continuous in space, discrete in time) for  curvature-type equations. See \cite{S1977} for an analogous game in combinatorics. 

The proof is standard under the framework of viscosity solution. For the reader's convenience, we provide an outline of the proof. 

\medskip

{\bf Step 1:} Define
\begin{eqnarray*}
\overline{u}(x,t):=\limsup_{y\rightarrow x,\ \tau\to 0,\ N\tau^2\to t}u_{\tau}(y,N\tau^2)\\
\underline{u}(x,t):=\liminf_{y\rightarrow x,\ \tau\to 0,\ N\tau^2\to t}u_\tau(y,N\tau^2).
\end{eqnarray*}
Then $\overline u$ is upper semicontinuous and $\underline {u}$ is lower semicontinuous. 
\begin{lem}\label{lem:boundarycontrol} For $t\in [0,1]$, 
$$
\min_{|y-x|\leq C\sqrt{t}}g(y)\leq \underline u(x,t)\leq \overline u(x,t)\leq \max_{|y-x|\leq C\sqrt{t}}g(y).
$$
Here $C$ is a constant depending only on $d$ and $\max_{\Rset^2}|V|$. 
\end{lem}
Proof: Fix $t>0$. Consider the game starting from $x_0\in \Rset^2$ with $(\tau, N)$, $N\tau^2\in [{t\over 2}, 2t]$. 

On the one hand, if player I adopts the following strategy similar to the exit strategy in \cite{KS1}: at each step $n$, player I chooses the direction $\eta_n$ such that 
$$
\eta_n\cdot (x_n-x_0)=0,
$$
then regardless of how player II plays, if $|x_n-x_0|\geq \sqrt{t}$, then
$$
|q+w|-|q|={2w\cdot q+|w|^2\over |q+w|+|q|}\ \Rightarrow\ |x_{n+1}-x_0|\leq |x_n-x_0|+ {C\tau^2}\left(1+{1\over \sqrt{t}}\right).
$$
Hence  
$$
|x_N-x_0|\leq C\sqrt{t} \quad \text{for $t\in [0,1].$}
$$
Then by the definition of $\overline u(x,t)$, 
$$ 
\overline u(x_0,t)\leq \max_{|y-x_0|\leq C\sqrt{t}}g(y).
$$

On the other hand, if player II employs the strategy as: at each step $n$, player II chooses the sign $b_n$ such that 
$$
b_n\eta_n\cdot (x_n-x_0)\leq 0
$$
after player I picks the direction $\eta_n$, then regardless what kind of strategy player I uses, similarly, we have that 
$$
|x_N-x_0|\leq C\sqrt{t}.
$$
This implies that 
$$
\underline u(x_0,t)\geq \min_{|y-x_0|\leq C\sqrt{t}}g(y).
$$
\qed
Since $g\in UC(\Rset^2)$, 
$$
\lim_{y\to x, t\to 0}\overline u(y,t)=\lim_{y\to x, t\to 0}\underline u(y,t)=g(x) \quad \text{uniformly for $x\in \Rset^2$}. 
$$

\medskip 

{\bf Step 2:} As in \cite{KS1}, we have that 
\begin{lem}
$\overline u(x,t)$ is a viscosity subsolution of 
$$
\overline u_t+\underline F(D^2\overline u, D\overline u)+V(x)\cdot D\overline u=0 \quad \text{on $\Rset^2\times (0, \infty)$}
$$
and
$\underline u(x,t)$ is a viscosity supersolution of 
$$
{\underline u}_t+\overline F(D^2\underline u, D\underline u )+V(x)\cdot D\underline u=0 \quad \text{on $\Rset^2\times (0, \infty)$}.
$$
\end{lem}
\qed

\medskip

{\bf Step 3:} Finally, the comparison principle Theorem \ref{theo:comparison} implies that $\underline u\geq \overline u$. Therefore, $\underline u= \overline u$ and (\ref{limit}) holds. \qed

\begin{rmk}\label{rmk:boundary} By the uniqueness of solution and Lemma \ref{lem:boundarycontrol}, if $n=2$ and $G(x,t)$ is the unique viscosity to equation (\ref{ge1}) subject to $G(x,0)=g(x)$,  Then for $t\in [0,1]$
$$
\min_{|y-x|\leq C\sqrt{t}}g(y)\leq G(x,t)\leq \max_{|y-x|\leq C\sqrt{t}}g(y).
$$
Of course, this can also be proved by pure PDE methods through comparison principle. 
\end{rmk}

\begin{rmk}\label{rmk:stochastic} Although the game is formulated in a deterministic way, 
it has some intrinsic stochastic features due to the different scales $\tau$ and $\tau^2$ 
so that the game trajectory is practically hard to analyze. 
Having the positive part (equivalently, allowing player I to choose $\eta_n=0$) renders 
the game more deterministic. As mentioned in \cite{KS2}, one of the two main questions 
about the game is whether it can be used to prove new results about PDE, 
which requires finding the right structure of PDE to work with. 
The curvature G-equation with integrable ambient fluid flows 
turns out to be a good candidate. The game formulation does allow us to 
make the best use of the underlying structure that is not naturally accessible 
by the relatively rough PDE approaches.  
See \cite{Liu2013} and \cite{HL} for other works related to applications of the game. 
\end{rmk}

\subsection{Stationary equation and reachability} 
As in \cite{KS1}, we may also consider stationary equations
\be\label{ge2}
\left(1-d\, \, \mathrm{div}\left({Du\over |Du|}\right)\right)_{+}|Du|+V(x)\cdot Du=\alpha.
\ee
for a constant $\alpha\geq 0$. 

\begin{defin} $u\in USC(\Omega)$ is called a viscosity subsolution of equation (\ref{ge2}) provided that given $\phi(x)\in C^2(\Rset^n)$, if for $x_0\in \Omega$
$$
u(x_0)-\phi(x_0)=\max_{x\in \Omega} (u(x)-\phi(x)),
$$
then
$$
\underline{F}\left(D^2\phi(x_0), D\phi (x_0)\right)+V(x_0)\cdot D\phi(x_0)\leq \alpha.
$$
\end{defin}

\begin{defin} $u\in LSC(\Omega)$ is called a viscosity supersolution of equation (\ref{ge2}) provided that given $\phi(x)\in C^2(\Rset^n)$, if for $x_0\in \Omega$
$$
u(x_0)-\phi(x_0)=\min_{x\in \Omega} (u(x)-\phi(x)),
$$
then
$$
\overline{F}\left(D^2\phi(x_0), D\phi (x_0)\right)+V(x_0)\cdot D\phi(x_0)\geq \alpha.
$$
\end{defin}

$u\in C(\Omega)$ is called a viscosity solution of (\ref{ge2}) if it is both a viscosity subsolution and a viscosity supersolution. Unlike the Cauchy problem, the stationary problem might not have a solution with a prescribed boundary data. 

Similar to the proof of Theorem \ref{theo:comparison}, we can prove the following comparison principle.

\begin{theo}\label{theo:comparison2} Assume that $\alpha>0$ and $\Omega$ is a bounded open set. Suppose that $u_1\in USC(\Omega)\cap L^{\infty}(\Omega)$ is a viscosity subsolution of equation (\ref{ge2}) on $\Omega$ and $u_2\in LSC(\Omega)\cap L^{\infty}(\Omega)$ is a viscosity supersolution of equation (\ref{ge2}) on $\Omega$. If for every point $x\in \partial \Omega$, 
$$
\limsup_{y\in \Omega\to x}u_1(y)\leq \liminf_{y\in \Omega\to x}u_2(y), 
$$
then
$$
u_1(x)\leq u_2(x) \quad \text{for $x\in \Omega$}. 
$$
\end{theo}

\begin{lem}\label{lem:mini} Suppose $\alpha>0$ and $u\in C(\overline \Omega)$ is a viscosity supersolution of equation (\ref{ge2}), then
$$
\min_{x\in \overline \Omega}u(x)=\min_{x\in \partial \Omega}u(x).
$$
\end{lem}

Proof: By the definition of viscosity supersolution, $u$ cannot attain minimum at $x_0\in \Omega$. Otherwise, we can use the constant function $\phi\equiv \min_{\bar \Omega}u$ as the test function and obtain
$$
0=\overline F(D^20, D0)+V(x_0)\cdot 0\geq \alpha>0, 
$$
which is absurd. Hence our conclusion holds. \qed

\begin{defin}\label{def:reach} Consider the game dynamics introduced in section 2.2. 
Let $S$ be subset of $\Rset^2$. We say that $S$ is reachable 
from $x$ (or $x$ can reach $S$) within time $T$ if for every 
open set $U$ satisfying $S\subset U$, there exist a sequence of 
positive numbers $\{\tau_m\}_{m\geq 1}$ such that $\lim_{m\to +\infty}\tau_m=0$ and 
for each fixed $\tau_m$, player I has a $U$-oriented strategy $\Gamma_{m1}$ such that 
regardless of how player II chooses his strategy $\Gamma_2$, under strategy $\Gamma_{m1}$, 
player I can force the associated game trajectory starting from $x$ to 
enter $U$ at $N(m, \Gamma_{m1}, \Gamma_2)$-th step  (i.e., $x_{N(m, \Gamma_{m1}, \Gamma_2)}\in U$), 
where  $N(m, \Gamma_{m1}, \Gamma_2)$ is a positive integer 
depending on $\tau_m,  \Gamma_{m1}, \Gamma_2$ and 
satisfying $N(m, \Gamma_{m1}, \Gamma_2)\tau_m^2\leq T$.

\end{defin}

\begin{rmk}\label{rmk:opentransition} To simplify notations, we usually omit the dependence of the game trajectory (in particular, the terminal point $x_{N}$) on $\tau$, $N$ and the strategy of player II. Moreover, by the above definition, it is clear that if a point $x$ can reach an open set $U$ within time $t_1$ and every point on $U$ can reach a set $S$ within $t_2$. Then the point $x$ can reach $S$ within time $t_1+t_2$. 
\end{rmk}

Hereafter,  for any set $E\subset \Rset^n$ and $t\geq 0$,  we denote by $E_t$  the image of $E$ under the $-V$ flow $\dot \xi=-V(\xi)$ at time $t$, i.e.,
$$
E_t=\{\xi_x(t)|x\in E, \  \dot \xi_x(s)=-V(\xi_x(s)),\  \xi_x(0)=x\}.
$$

\begin{lem}\label{lem:reachcontrol} Suppose that $u\in UC(\Rset^2)$ is a 
viscosity subsolution of equation (\ref{ge2}) on $\Rset^2$. 
For $x_0\in \Rset^2$, if there exists $T_0>0$ such that a bounded set $S$ is reachable from $x_0$ within time $T_0$, then 
$$
u(x_0)\leq \max_{y\in \cup_{t\in [0,T_0]}\bar {S_t}}u(y)+\alpha T_0.
$$
In particular,  if the set $S$ is $-V$  flow invariant (i.e., $S_t=S$ for all $t\geq 0$), then
$$
u(x_0)\leq \max_{y\in \bar S}u(y)+\alpha T_0.
$$
\end{lem}
Proof: let $w(x,t)=u(x)-\alpha t$. Then $w(x,t)$ is a viscosity subsolution to equation (\ref{ge1}) with initial data $w(x,0)=u(x)$. Now let $G(x,t)\in C(\Rset^2\times [0, \infty))$ be the unique viscosity solution to equation (\ref{ge1}) with initial data $G(x,0)=u(x)$. Owing to the comparison principle Theorem \ref{theo:comparison},
$$
w(x,t)\leq G(x,t).
$$
For $r>0$ and $E\subset \Rset^2$, let $D_r(E)=\{x\in \Rset^2|\ d(x, E)<r\}$.  Fix $r>0$.  Due to the game formulation of $G(x,t)$ and the Definition \ref{def:reach} of reachability, there are a sequence of positive numbers $\{\tau_m\}_{m\geq 1}$ and a sequence of positive integers $\{N_m\}_{m\geq 1}$ such that  $\lim_{m\to +\infty}\tau_m=0$ and  for each $\tau_m$,  player I has a $D_r(S)$-oriented strategy to drive the game trajectory to enter $D_r(S)$ at the $N_m$-th step for $N_m\tau_m^2\leq T_0$  regardless of how player II plays.  Here for convenience, we omit the dependence of $N_m$ on the concrete strategy of player II since what really matters here is the upper bound $N_m\tau_m^2\leq T_0$. 

\medskip

Next,  starting from $x_{N_m}\in D_r(S)$, player I chooses $\eta=0$ for $k$ more steps. Here $k$ is the 
first whole number such that $T_0\geq (N_m+k)\tau_m^{2}\geq T_0-\tau_m^2$.  

Let $J_m$ be the integer part of $T_0\over \tau_m^2$. Then the above argument says that player I has a strategy such that regardless of how player II plays,  player I can steer the game trajectory into
$$
D_{C\tau_{m}^2}(E_r).
$$
at the $J_m$-th step for $E_r=\cup_{0\leq t\leq T_0}(D_r(S))_t$.

Hence for the value function defined in (\ref{gamedef}), due to $G(x,0)=u(x)$, 
$$
u_{\tau_m}(x_0, J_m\tau_{m}^{2})\leq \max_{y\in \overline {D_{C\tau_m^2}(E_r)}}u(y)
$$
$$
G(x_0,T_0)=\lim_{m\to +\infty}u_{\tau_m}(x_0, J_{m}\tau_{m}^{2})\leq  \max_{y\in \overline E_r}u(y).
$$
So
$$
u(x_0)-T_0\alpha=w(x_0,T_0)\leq G(x_0,T_0)\leq \max_{y\in \overline E_r}u(y).
$$
Accordingly, our conclusion follows by sending $r\to 0$. \qed

\begin{lem}\label{lem:flowcontrol}  Suppose that $u\in UC(\Rset^n)$ is a viscosity subsolution of equation (\ref{ge2}) on $\Rset^n$.  Let $\xi: [0,\infty)\to \Rset^n$ satisfy $\dot \xi(s)=-V(\xi(s))$.  Then for $t_1<t_2$,
$$
u(\xi(t_1))\leq u(\xi(t_2))+\alpha (t_2-t_1).
$$
\end{lem}

Proof: This follows directly from the two facts below:

\medskip

(1)  $w(x,t)=u-\alpha t$ is also a viscosity subsolution of the  transport equation
$$
w_t+V(x)\cdot Dw=0 \quad \text{subject to $w(x,0)=u(x)$},
$$

(2)  the unique solution to
$$
\begin{cases}
v_t+V(x)\cdot Dv=0\\[3mm]
v(x,0)=u(x)
\end{cases}
$$ 
is given by the representation formula: $v(x,t)=u(\xi_x(t))$ for $\dot \xi_x(s)=-V(\xi_x(s))$ subject to $\xi_x(0)=x$. \qed

\section{Proof of Theorem \ref{theo:main}}

Let $p\in \Rset^n$ be a fixed unit vector. By the standard Perron's method \cite{CL1992}, for any given $\lambda>0$, the following equation has a unique continuous $\Zset^n$-periodic viscosity solution $v=v_{\lambda}\in C(\Rset^n)$. 
$$
\lambda v+ \left(1-d\Div{\frac{p+Dv}{|p+Dv|}}\right)_{+}|p+Dv|+V(x)\cdot (p+Dv) =0 \quad \text{in $\Rset^n$}.
$$
The above equation also has a comparison principle whose proof is similar to that of Theorem \ref{theo:comparison}. 
{To prove Theorem 1.1, our main
task will be to show} that when $n=2$ and $V=A(DH)^{\perp}$ is the 2D cellular flow, 
there exists a positive constant $\overline H(p)$ such that 
$$
\lim_{\lambda\to 0}\lambda v_{\lambda}(x)=-\overline H(p) \quad \text{uniformly on $\Rset^2$}. 
$$
Here we omit the dependence of $\overline H$ on the flow intensity $A$. 

Throughout this section, $C$ represents a constant depending only on $d$ and $V$. By maximum principle, it is easy to see that 
\be\label{v-max}
\max_{x\in \Rset^2}|\lambda v_\lambda(x)|\leq 1+\max_{\Rset^2}|V|.
\ee

{ We believe that the equi-continuity of  the family of functions $\{v_{\lambda}\}_{\lambda>0}$  does not hold  due to the degeneracy of  the curvature term and the lack of coercivity.}  Instead, we will show that 
$$
\max_{x,y\in \Rset^2}|v_\lambda(x)-v_{\lambda}(y)|\leq C.
$$
Our strategy is (1) establishing partial reachability from the associated game dynamics based on the special structure of the cellular flow to apply Lemma \ref{lem:reachcontrol} and Lemma \ref{lem:flowcontrol}, then (2) using the minimum value principle Lemma \ref{lem:mini} to compensate for the lack of full reachability. The proof can be viewed as a combination of Lagrangian and Eulerian approaches.  The  game trajectory under player I's strategy more or less mimics the reverse of the propagation route of flame.

By adopting the proof of Theorem 4 in \cite{AB}, we {first establish the connection between $G(x,t)$ and $v_\lambda$.  Combining with partial reachability, this will lead to a negative upper bound of $\lambda v_{\lambda}(x)$ for small $\lambda$ in Corollary 3.2 later, which allows us to apply the minimum value principle.} 

\begin{lem}\label{lem:ergodic-version} Let $G(x,t)$ be the unique solution of (\ref{ge1}) with $G(x,0)=p\cdot x$. Suppose that there exists $\beta>0$ such that for all $(x,t)\in \Rset^n\times [0, \infty)$
$$
G(x,t)-p\cdot x\leq -\beta t +C.
$$
Then 
$$
\max_{x\in \Rset^n}\lambda v_{\lambda}(x)< -{\beta \over 2}+ \lambda C
$$
Here $C$ represents a constant depending only on $d$ and $V$. 
\end{lem}
Proof: By comparing $G(x,t)$ and $p\cdot x\pm Mt$ for a suitable constant $M$ depending only on $d$ and $V$, the comparison principle Theorem \ref{theo:comparison} implies that 
$$
|G(x,t)-p\cdot x|\leq Mt \quad \text{for all $(x,t)\in \Rset^n\times [0, \infty)$}. 
$$
Apparently, $\beta\leq M$. 

In addition, by periodicity of $V$, $G(x+\vec{l}, t)-p\cdot \vec{l}$ is also a solution of (\ref{ge1}) with intial data $p\cdot x$ for any $\vec{l}\in \Zset^n$. Then uniqueness implies that $G(x+\vec{l}, t)-p\cdot \vec{l}=G(x,t)$, equivalently, $G(x,t)-p\cdot x$ is $\Zset^n$-periodic for $x$. By the assumption, we may choose $T_0>0$ such that 
$$
G(x,t)-p\cdot x\leq -{2\beta t\over 3} \quad \text{for $(x,t)\in \Rset^n\times [T_0, \infty)$}.
$$

{\bf Step 1:} Choose $f(t)\in C^{\infty}([0, \infty))$ such that 
$$
f'(t)\leq -{\beta\over 2}
$$
and
$$
f(t)=
\begin{cases}
-2Mt \quad \text{when $t\geq 2T_0$}\\[3mm]
-{\beta t\over 2} \quad \text{when $t\in [0, T_0]$}.
\end{cases}
$$

Define
$$
h(x)=\min_{t>0} (G(x,t)-p\cdot x- f(t))>-\infty.
$$
Clearly, $h(x)<0$. Hence the minimum is attained for some $t_x\in (0, 2T_0]$. Therefore, $h(x)$ is periodic, continuous and is a viscosity supersolution of
$$
\left(1-d\Div{\frac{p+Dh}{|p+Dh|}}\right)_{+}|p+Dh|+V(x)\cdot (p+Dh) ={\beta\over 2} \quad \text{on $\Rset^n$}.
$$

\medskip

{\bf Step 2:} Then
$$
h_{\lambda}(x)=h(x)-\min_{x\in \Rset^n}h(x)-{\beta\over 2\lambda}
$$
is a viscosity supersolution of
$$
\lambda h_{\lambda}(x)+\left(1-d\, \Div{\frac{p+D h_\lambda}{|p+D h_\lambda|}}\right)_{+}|p+D h_{\lambda}|+V(x)\cdot (p+D h_{\lambda})=0.
$$
Accordingly, by comparison principle,
$$
v_{\lambda}(x)\leq h_{\lambda}(x),
$$
which implies that 
$$
\max_{x\in \Rset^n}\lambda v_{\lambda}(x)\leq -{\beta\over 2}+2\lambda\max_{ x\in \Rset^n}|h(x)|.
$$
Hence our conclusion holds. \qed

Hereafter, we let $x=(x_1,x_2)\in \Rset^2$ and 
$$
V(x)=A(DH)^{\perp}=A(-\cos(x_2)\sin (x_1), \ \cos(x_1)\sin (x_2)).
$$
for $A>0$. Keep in mind that our game is using
$$
-V(x)=A(\cos(x_2)\sin (x_1), \ -\cos(x_1)\sin (x_2)).
$$

Write
$$
Q=[0,\pi]\times [0,\pi], \quad Q_\mu=\{x\in Q| \ H(x)>\mu\}
$$
and
$$
\Gamma_\mu=\{x=(x_1,x_2)\in \Rset^2|\ \min\{|x_1|, |x_2|, |x_1-\pi|, |x_2-\pi|\}<\mu\}.
$$
for $\mu\in (0,1]$. Note that $Q_\mu$ is $-V$ flow invariant. 
\begin{center}
\includegraphics[scale=0.5]{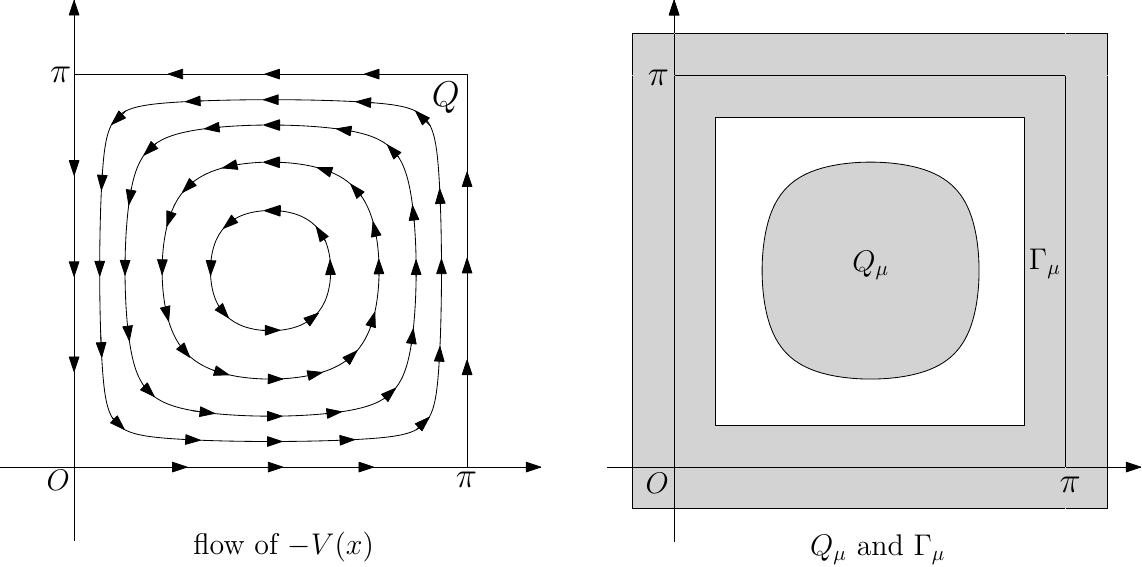} 
\captionof{figure}{Flow of $-V$, two domains $Q_{\mu}$ and $\Gamma_{\mu}$.}
\end{center}

We will utilize the special structure of the cellular flow to establish reachability. 

\begin{lem}\label{lem:interior-reach} Suppose $\mu \in [0,1)$.   For $P_1\in Q_\mu$,   the level curve $\{x\in Q|\ H(x)=\mu\}$  is reachable from $P_1$ within time $3\sqrt{2}$. In particular, for $\mu>0$,   every point $P_2\in \{x\in Q|\ H(x)=\mu\}$ is reachable from $P_1$ within time $\leq C(1+|\log \mu|)$ for a constant $C$ depending only on  $V$.

\end{lem}

Proof: First we prove the reachability to the level curve. Given the game step size parameter $\tau$, recall that the game dynamics is 
$$
\begin{cases}
X_{n+1}=X_n+\tau\sqrt{2d}b_n\eta_n+\tau^2\eta_{n}^{\perp}-\tau^2V(X_n)\\
X_0=P_1
\end{cases}
$$
for $n\geq 0$. Here $\eta_n\in \overline B_1(0)=\{v\in \Rset^2|\ |v|\leq 1\}$ and $b_n\in \{-1,1\}$. 

\medskip

Case 1: $P_1\not=({\pi\over 2}, {\pi\over 2})$. Then $|DH(P_1)|>0$. Player I chooses the strategy as follows:  Let $X_0=P_1$. 

\medskip

 At each step $n\geq 0$, if $H(X_n)>H(P_2)$, player I chooses 
$$
\eta_n={V(X_n)\over |V(X_n)|}.
$$
Then
$$
\eta_{n}^{\perp}=-{DH(X_n)\over |DH(X_n)|}.
$$
Then regardless of how player II moves, we have that 
$$
\begin{aligned}
&H(X_{n+1})-H(X_n)=H\left(X_n+\tau b_n\sqrt{2d}\eta_n+\tau^2\eta_n^\perp-\tau^2V(X_n)\right)-H(X_n)\\
=&\inner{DH}{\tau b_n\sqrt{2d}\eta_n+\tau^2\eta_n^\perp-\tau^2V(X_n)}+\frac{1}{2}\inner{D^2H\tau b_n\sqrt{2d}\eta_n}{\tau b_n\sqrt{2d}\eta_n}+O(\tau^3)\\
=&\tau^2\left[-\norm{DH}+d\inner{D^2H\Dpn{H}}{\Dpn{H}}\right]+O(\tau^3).
\end{aligned}
$$
Since $H(X_n)=\sin{x_{1n}}\sin{x_{2n}}$ for $X_n=(x_{1n},x_{2n})$, we have
$$
-\max_{\Rset^2}||D^2H||\leq \inner{D^2H\Dpn{H}}{\Dpn{H}}=\frac{ \sin (x_{1n}) \sin (x_{2n}) (\cos (2 x_{1n})+\cos (2 x_{2n})+2)}{\cos (2x_{1n})\cos(2x_{2n})-1}\leq 0.
$$
Here the negative sign is essentially due to the convexity of the level curve of $H$ instead of the specific form of $H$. So 
$$
H(X_{n+1})-H(X_n)\leq -|DH(X_n)|\tau^2+O(\tau^3).
$$
Let $a_n=H(X_n)$. Then by Lemma \ref{lem:gradient-control},
$$
a_{n+1}\leq a_n-\sqrt{2(a_n-a_n^{2})}\tau^2+O(\tau^3).
$$
This implies that the decreasing rate of $a_n=H(X_n)$ when $n$ increases is no slower than the 
decreasing rate along the ODE
$$
\dot s(t)=-\sqrt{2s(t)(1-s(t))} \quad \text{with $s(0)=H(P_1)\in (0,1)$},
$$
which has a unique solution when $s(t)\in (0,1)$.  {Assume that $s(t)>0$ for $t\in (0,t_0)$ and $s(t_0)=0$.  Clearly  $s(t)$ is strictly decreasing for $t\in [0,t_0)$. Since $\sqrt{2s(t)(1-s(t))}\geq \sqrt{1-s(t)}$ when  $s(t)\in ({1\over 2}, 1)$ and $\sqrt{2s(t)(1-s(t))}\geq \sqrt{s(t)}$ when  $s(t)\in (0,{1\over 2}]$, we have that 
$$
t_0\leq \int_{{1\over 2}}^{1}{1\over \sqrt{1-s}}\,ds+\int_{0}^{1\over 2}{1\over \sqrt{s}}\,ds=2\sqrt{2}. 
$$
In particular,  the upper bound is independent of the value of $H(P_1)$.}  
Hence a simple calculation shows  that when $\tau$ is small enough, after at most $n_1\leq {3\sqrt{2}\over \tau^2}$ steps, 
we have that 
$$
H(X_{n_1})\leq \mu
$$
for the first time. In particular, $d(X_{n_1}, \{x\in Q|\ H=\mu\})\leq C\tau$. Hence by Defintion \ref{def:reach}, the needed time to reach the level curve $T_1\leq 3\sqrt{2}$. 

\begin{center}
\includegraphics[scale=0.3]{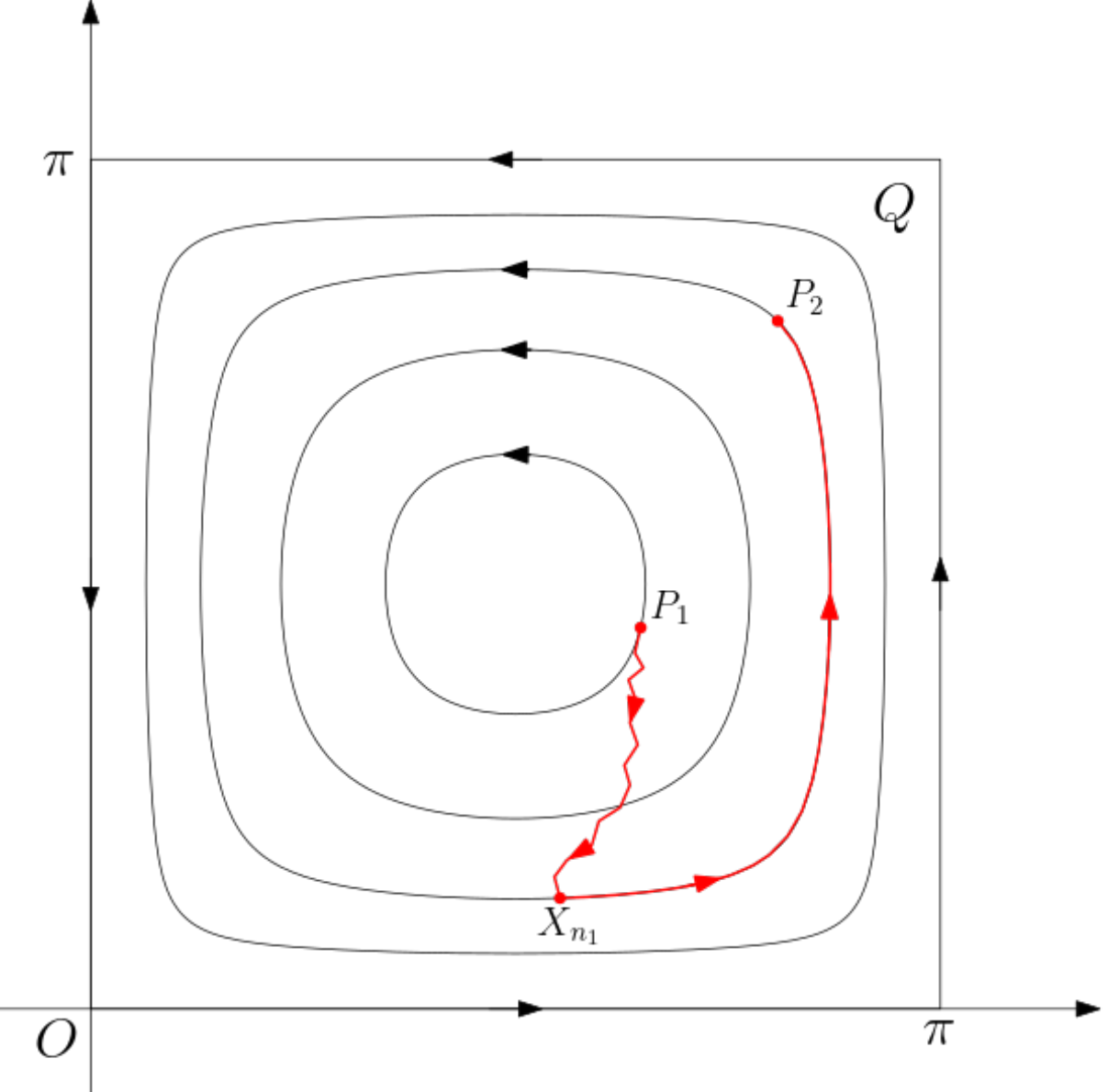} 
\captionof{figure}{Game trajectory of Lemma \ref{lem:interior-reach}.}
\end{center}

Case 2: $P_1=({\pi\over 2}, {\pi\over 2})$. First travel a little bit away from $P_1$.   Choose $r>0$ such that $\max_{\overline B_r(P_1)}|V|\leq {1\over 4}$. Then player I can use the exit strategy as in the proof of Lemma \ref{lem:boundarycontrol}. It is easy to see that player I could reach $\partial B_r(P_1)$ within time $Cr^2$. Then it goes to Case 1. 

\medskip  

 Next we show the reachability to every point $P_2$ on the level curve $\{x\in Q|\ H(x)=\mu\}$. After the game trajectory reaches the level curve , player I can just choose $\eta=0$, i.e.,  follow the $-V$ flow to travel along the level curve  to reach every point $P_2$ on the curve.  When $\mu$ is close to 1, the traveling time around the level curve  is near $2\pi$. When $\mu$ is close to 0, the traveling time is bounded by $O(|\log\mu|)$.  Hence the conclusion holds. 
\qed

\begin{rmk} The reachability established in the above lemma is only one way: $P_1$ is NOT reachable from $P_2$  if $H(P_1)$ is close to 1 since the curvature on the level curve will become very big and surpass 1. This is different from the invicid case ($d=0$) where two points are mutually reachable. 
\end{rmk}

\begin{lem}\label{lem:gradient-control} For $x=(x_1,x_2)\in [0, \pi]\times [0, \pi]$,
$$
|DH(x)|\geq \sqrt{2(H(x)-H^2(x))}. 
$$
\end{lem}

Proof: A direction computation shows that
$$
\begin{aligned}
|DH(X)|&=\sqrt{\sin^2 x_1+\sin^2 x_2-2\sin^2 x_1\sin^2 x_2}\\
&\geq \sqrt{2\sin x_1\sin x_2-2\sin^2 x_1\sin^2 x_2}\\
&=\sqrt{2(H(x)-H^2(x))}.
\end{aligned}
$$
\qed

{ We would like to point out that the convexity of level sets of $H$ and the above inequality provide technical convenience in the proof, but they are not essential in obtaining the existence of $\overline H(p)$. See Section 4 for more details. }

\begin{lem}\label{lem:boundary-reach}
There exists $\mu_0>0$ and $T_{0}>0$ depending only on $d$ and $V$ such that the set $Q_{2\mu_0}$ is reachable from every point $x\in \Gamma_{\mu_0}$ within time $T_{0}$. 
\end{lem}

Proof: It is enough to prove the above conclusion when $\Gamma_{\mu_0}$ is replaced by the set
$$
\mathcal{A}_{\mu_0}=\left[-\mu_0, {2\pi\over 3}\right]\times [-\mu_0,\mu_0] 
$$
since the proof for other pieces are similar and we can just choose the smallest $\mu_0$.

\medskip

{\bf Step 1:} We first consider points near the corner $O=(0,0)$.  For $\alpha>0$, set
$$
\Omega_\alpha=(\alpha, 1)\times \left(-{1\over 2}, \ {1\over 2}\right ).
$$
{{Note that  $V(x)\cdot (1,0)=0$ for all $x\in \{0\}\times \Rset$.}} By applying Lemma \ref{lem:perreach} for $S=(0, 1)\times \left(-{1\over 2}, \ {1\over 2}\right )$ and $\Omega=\Omega_\alpha$, 
we have that there is $\alpha_1>0$ such that $\Omega_{2\alpha_1}$ is reachable from every point $x\in [-\alpha_1,\alpha_1]^2$ within time 1.

\medskip

{\bf Step 2:} Next we look at points on the line segment
$$
L_1=\left\{(x_1,0)|\ x_1\in \left[{\alpha_1\over 2}, {2\pi\over 3}\right]\right\}.
$$
By applying Lemma \ref{lem:perreach} for $S=Q$ and $\Omega=Q_{\mu}$, we deduce that there is $\mu_0\in (0, {\alpha_1\over 2})$  such that $Q_{2\mu_{0}}$ is reachable by every $x\in L_{1,\mu_0}= L_1\times [-\mu_0,\mu_0]$ within time 1.

\medskip

{\bf Step 3:} Next, given a point $P_1\in [-\alpha_1,\alpha_1]^2$, by step 1, player I has a 
strategy to push the game trajectory to some point $P_2\in \Omega_{2\alpha_1}$ within time $1$. 
If $P_2\notin L_{1,\mu_0}$, then $|H(P_2)|\geq \sin \mu_0\sin (2\alpha_1)$. Hence, 
by Lemma \ref{lem:interior-reach}, $P_2$ can reach $L_{1,\mu_0}$ within time $t_1$. 
Thus by step 2, $P_1$ can reach $Q_{2\mu_{0}}$ within time $1+1+t_1=2+t_1$. 
\begin{center}
\includegraphics[scale=0.7]{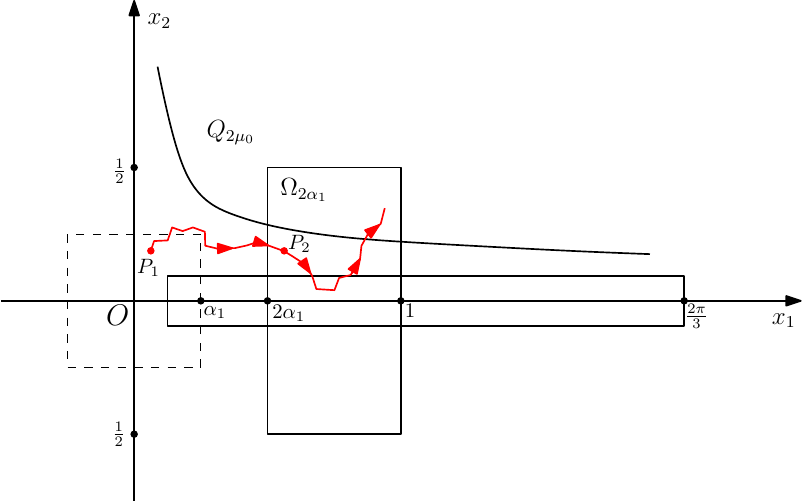} 
\captionof{figure}{Game trajectory of Lemma \ref{lem:boundary-reach}.}
\end{center}

Accordingly, every point on $\mathcal{A}_{\mu_0}\subset \left([-\alpha_1,\alpha_1]^2\cup L_{1,\mu_0}\right)$ can reach $Q_{2\mu_{0}}$ within time $2+t_1$. Note that $\mu_0$ (hence $t_1$ as well) depends only on $d$ and $V$. 

\qed

Given a point $X_0=(a_0, b_0)$ and $\tilde \delta>0$, let $\gamma(s)\in C([0,1], \Rset^2)$ be a continuous simple curve satisfying
$$
\gamma((0,1))\in (a_0-\tilde \delta, a_0+\tilde \delta)\times (-\infty, b_0)
$$
and $\gamma(0)\cdot (1,0)=a_0-\tilde \delta$, $\gamma(1)\cdot (1,0)=a_0+\tilde \delta$. Write
$$
J_r=\left(a_0-{\tilde \delta\over 2}, \ a_0+{\tilde \delta\over 2}\right)\times \left(b_0+r, \ b_0+r+r'\right)
$$
for some $r,r'>0$. Note that the curve $\gamma$ divides the strip $[a_0-\tilde \delta, a_0+\tilde \delta]\times \Rset$ into two connected components $I_1$ and $I_2$. Let $I_1$ be the lower component.

\begin{lem}\label{lem:cross} Given a point $x=(x_1,x_2)\in I_1$, if $x$ can reach $J_r$ within $t_0\leq {\tilde \delta\over 3M}$, then $x$ can also reach $\gamma([0,1])$ within $t_0$. Here $M=1+\max_{\Rset^2}|V|$. 
\end{lem}

Proof: According to Definition \ref{def:reach} of the reachability, player I has a $J_r$-oriented strategy to move the game trajectory into $J_r$ within time $t_0$ regardless of how player II responses. Let $U$ be an open set such that $\gamma([0,1])\subset U$ and
$$
\nu=d(\gamma([0,1]), \partial U)<1.
$$

We claim that the corresponding game trajectory of player I's $J_r$-oriented strategy with a game step size parameter $\tau<{\nu\over M+d}$ must enter $U$ before it arrives at $J_r$. If not, then at some moment $t_1<t_0$, the game trajectory arrives at 
$$
\{x=(x_1,x_2)|\ |x_1-a_0|>\tilde \delta\}.
$$
Without loss generality, we assume that it arrives at $\{x_1<a_0-\tilde \delta\}$. 
\begin{center}
\includegraphics[scale=0.4]{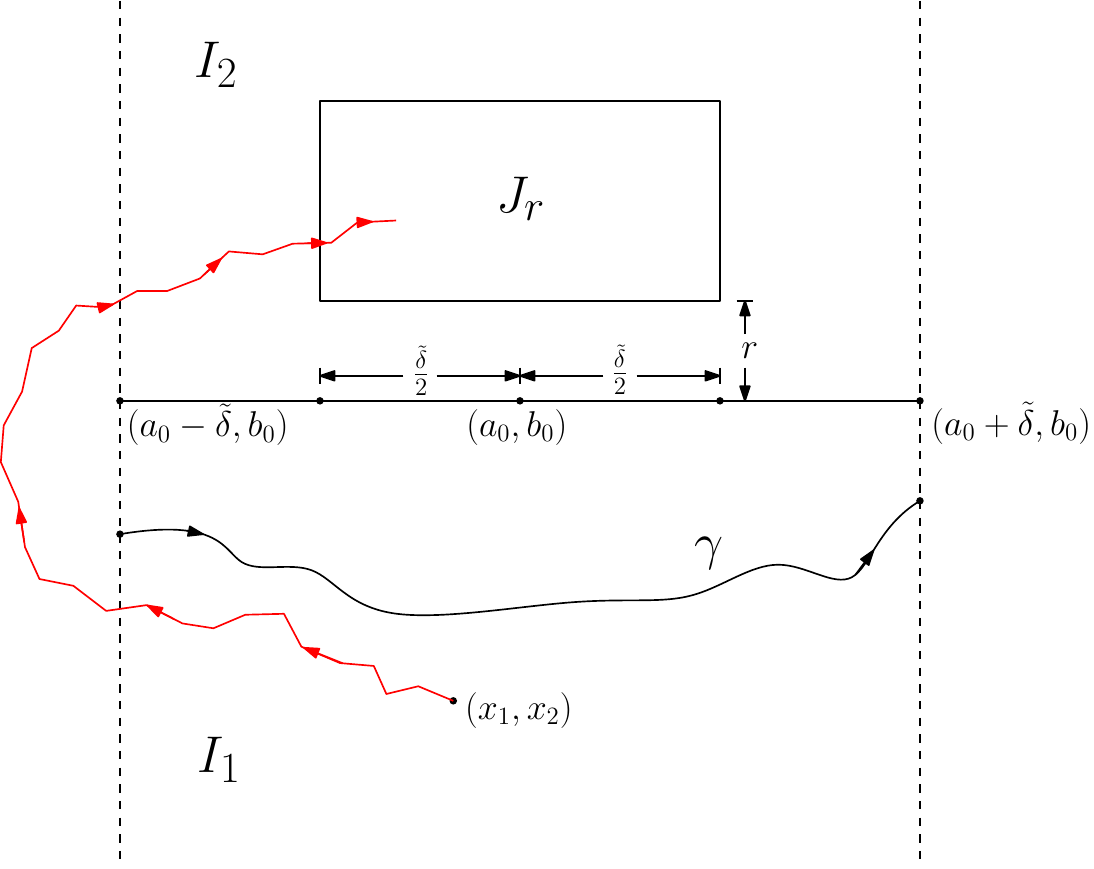} 
\captionof{figure}{Game trajectory of Lemma \ref{lem:cross}.}
\end{center}
Then player II can simply employ the strategy: at each step $n$, after play I picks the direction $\eta_n$, player II chooses the sign $b_n$ such that 
$$
b_n\eta_n\cdot (1,0)\leq 0.
$$
Then the only remaining force that can push the game trajectory to move along the positive horizontal direction is the mean burning speed 1 and the flow $-V$. Therefore, it will take the game at least $\tilde \delta\over 2M$ time to reach $J_r$. Since $t_0 \leq {\tilde \delta\over 3M}< {\tilde \delta\over 2M}$, this is a contradiction. \qed

\begin{cor}\label{cor:cross-reach}For any $0<\theta<{\pi\over 8}$, every point in the set
$$
Z_\theta=\left((\theta, \pi-\theta)\times \left[0, \pi\right]\right)\cup \left(\left[0, \pi\right]\times (\theta, \pi-\theta)\right)
$$
is reachable from the center $P_0=({\pi\over 2}, {\pi\over 2})$ within time $t_{\theta}$ depending only on $\theta$, $V$ and $d$. 
\end{cor}

Proof:  Let $M=1+\max_{\Rset^2}|V|$.  Due to symmetric structure and Lemma \ref{lem:interior-reach}, it suffices to show that for some $\delta_\theta>0$, every point in the set
$$
S_{\theta}=\left([\theta, \pi-\theta]\times \left[0,{\pi\over 2}\right]\right)\cap \{0\leq H< \delta_{\theta}\}
$$
is reachable from $P_0$ within time $t_{\theta}$ depending only on $\theta$, $V$ and $d$.

By applying Lemma \ref{lem:perreach} for $S=\left({3\theta\over 8},\ {5\theta\over 8}\right)\times \left(-{\theta\over 4},\ 0 \right)$ and
$$
\Omega_m=\left({3\theta\over 8},\ {5\theta\over 8}\right)\times \left( -{1\over m}-{\theta\over 4},\ -{1\over m} \right)
$$
we can show that there exist $m\in \Nset$ and $X_m=({\theta\over 2}, {1\over m})$ such that ${1\over m}<{\theta\over 4}$ and every point $x\in B_{1\over 4m}(X_m)$ can reach $\Omega_m$ within time $\theta\over 12M$. Now fix $m$ and set
$$
\delta_\theta={1\over 2}\min_{x\in \overline B_{1\over 4m}(X_m)}H(x).
$$
Applying Lemma \ref{lem:cross} with $\tilde \delta={\theta\over 4}$ to level curves of the stream function $H$ within the strip $\left[{\theta\over 4}, {3\theta\over 4}\right] \times \Rset$, we deduce that for any $\mu\in [0, \delta_\theta]$, every point on $B_{1\over 4m}(X_m)$ can reach the curve
$$
\xi_{\mu}=\{H=\mu\}\cap \left( \left[{\theta\over 4}, {3\theta\over 4}\right]\times \left(-\infty,\ {\pi\over 2}\right] \right)
$$
within time $\theta\over 12M$. 

Next we show that every point on $S_{\theta}$ is reachable from $P_0$. 

In fact, suppose that $P_1\in S_{\theta}$. Let $\mu_1=H(P_1)\in [0, \delta_\theta)$. 

Step 1: Thanks to Lemma \ref{lem:interior-reach}, $P_0$ can reach $X_m$ within time $t_1$. By the definition of reachability, the game trajectory can arrive at some point in $Y_1\in B_{1\over 4m}(X_m)$ within $t_1$.

\medskip

Step 2: By the above discussion and the definition of reachability, starting from $Y_1$, for any $k\geq 2\in \Nset$, the game trajectory can arrive at a point $Y_k$ such that $d(Y_k, \xi_{\mu_1})<{1\over k}$ within time $t_2={\theta\over 12M}$.

\medskip

Step 3: Starting from $Y_k$, at each step $n$, player I chooses $\eta_n=0$ (i.e., just following the flow $-V$). Then within time $t_3$, the game trajectory will arrive at some point as close to $P_1$ as we want when $k\to +\infty$. Note that $|V|$ has a positive lower bound depending on $\theta$ on $\left([{\theta\over 4}, \pi-\theta]\times \left[0,{\pi\over 2}\right]\right)$. So $t_3$ has a upper bound depending only on $\theta$ and $V$. 

\medskip

Hence $P_0$ can reach any neighborhood of $P_1$ within time $t_1+t_2+t_3$. Note that each $t_i$ depends only on $\theta$, $V$ and $d$. 

\begin{center}
\includegraphics[scale=0.4]{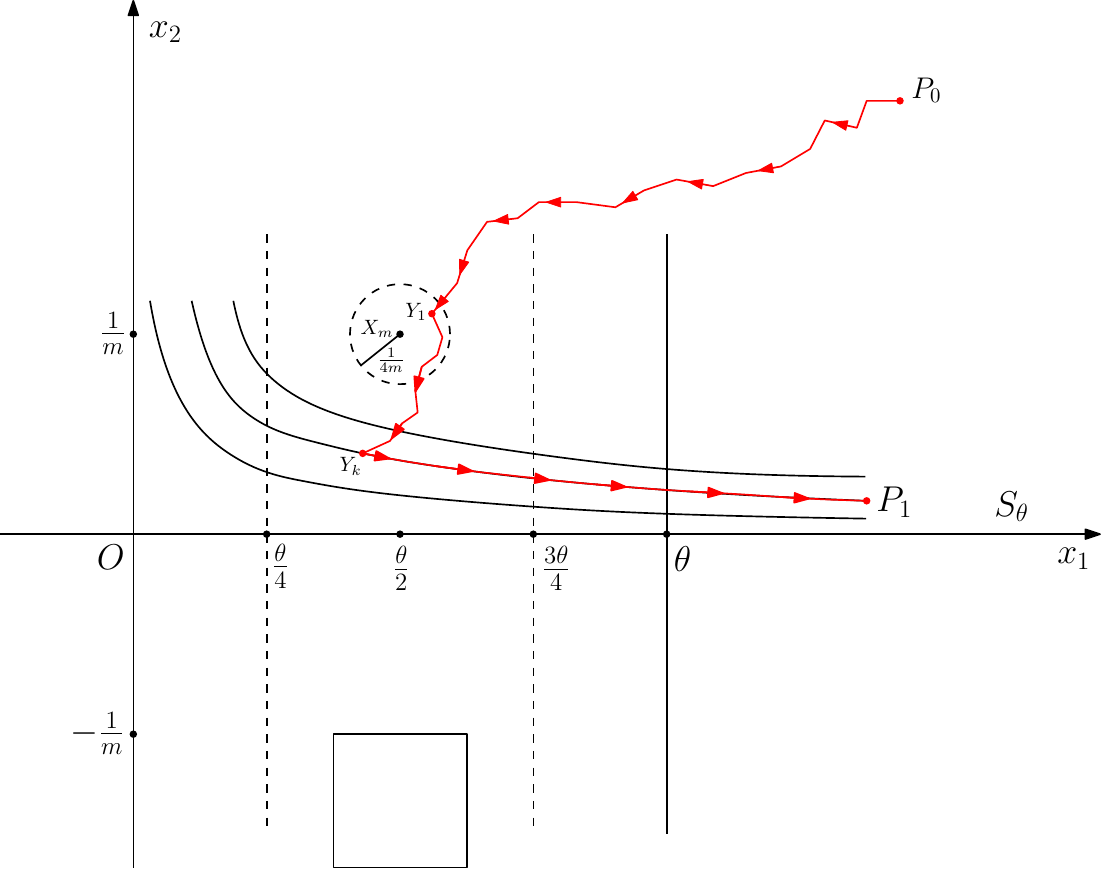} 
\captionof{figure}{Game trajectory of Corollary \ref{cor:cross-reach}.}
\end{center}
\qed

Recall that $Q=[0,\pi]^2$. Consider the interior of four cells in $[-\pi, \pi]^2$: $U_1=(0, \pi)\times (0,\pi)$, $U_2=U_1-(\pi,0)$, $U_3=U_1-(0,\pi)$ and $U_4=U_1-(\pi,\pi)$. Clearly, by similar proofs, the corresponding versions of previous results 
on $U_1$ could also be established for $U_i$ ($i=2,3,4$). 

Below is a transition property from one cell to another. 

\begin{lem}\label{lem:direction-reachability} There exists $\beta>0$ depending only on $d$ and $V$ such that every point $P\in U_1$ can reach the set $U_i$ within time $\beta$ for $i=2,3,4$.

\end{lem}

Proof: It suffices to prove this for $U_2$. The others are similar. Fix $P\in U_1$. 
Choose $\mu_0$ from Lemma \ref{lem:boundary-reach} that works for the corresponding statements for all four cells in $[-\pi,\pi]^2$. Choose the largest $\tilde \mu\in (0, \mu_0)$ such that
$$
\{x\in Q|\ H(x)\leq \tilde \mu\}\subset \overline{\Gamma}_{\mu_0}.
$$

Case 1: If $H(P)\geq \tilde \mu$, by Lemma \ref{lem:interior-reach}, $P$ can reach at some point 
$$
P_1\in J= \left\{x\in Q|\ {\tilde \mu\over 2}<H(x)<{\tilde \mu}\right\}
$$
within time $t_1$. Then by following the flow $-V$, $P_1$ can reach a point $P_2\in J\cap ((0,\mu_0)\times (0,\pi))$ within time $t_2$. Since 
$$
(0,\mu_0)\times (0,\pi)\subset \Gamma_{\mu_0}-(\pi,0), 
$$
$P_2$ can reach $U_2$ within time $T_0$ by applying Lemma \ref{lem:boundary-reach} to $U_2$.

\medskip

Case 2: If $H(P)<\tilde \mu$, then $P\in \Gamma_{\mu_0}$. Owing to Lemma \ref{lem:boundary-reach}, $P$ can arrive at some $\tilde P_1\in Q_{2\mu_0}$ within time $t_3$ , which goes back to case 1. 

In the above, $t_1$, $t_2$ and $t_3$ depend only on $d$ and $V$. Let $\beta=t_1+t_2+t_3+T_0$. \qed
\begin{center}
\includegraphics[scale=0.4]{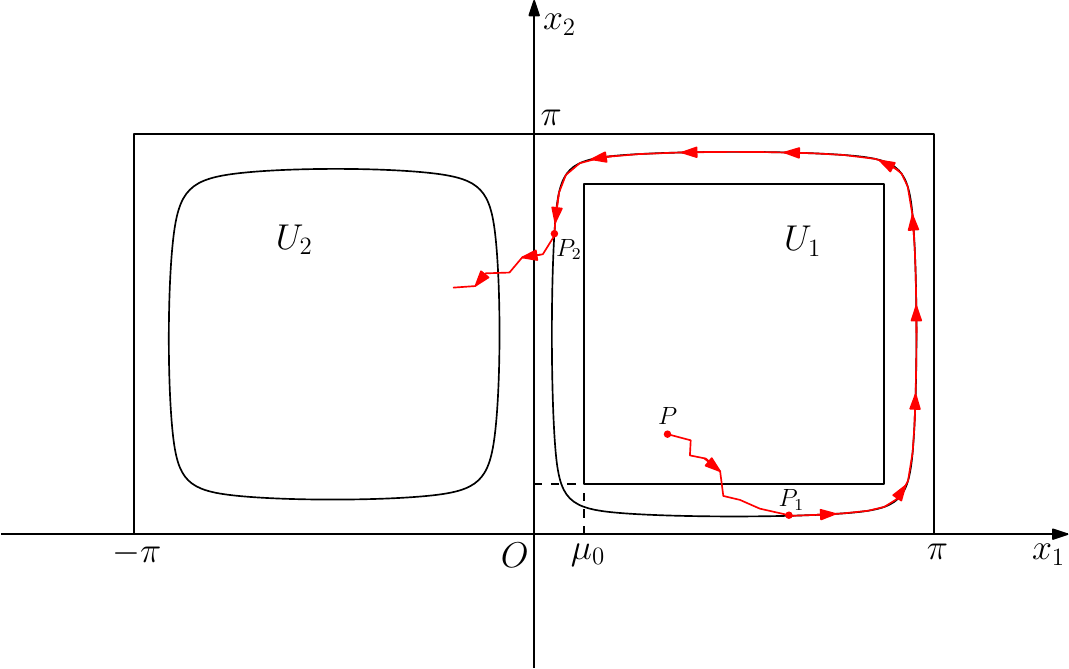}
\captionof{figure}{Case 1 game trajectory of Lemma \ref{lem:direction-reachability}. }
\end{center}
\begin{center}
\includegraphics[scale=0.27]{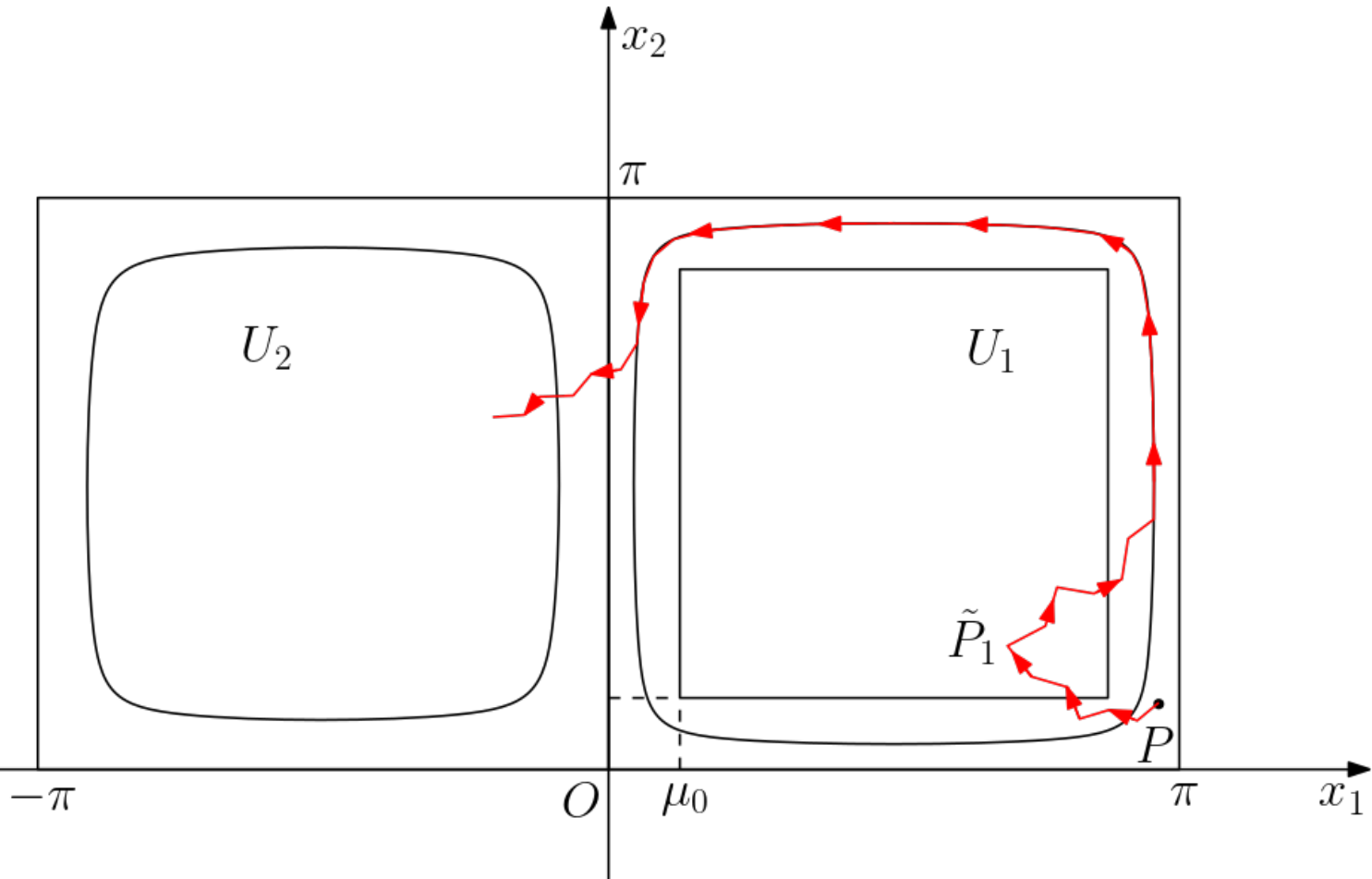} 
\captionof{figure}{Case 2 game trajectory of Lemma \ref{lem:direction-reachability}. }
\end{center}

\qed

\begin{lem}\label{lem:cauchy-limit} Let $p\in \Rset^2$ be a unit vector and $G(x,t)$ be the unique solution of (\ref{ge1}) with $G(x,0)=p\cdot x$. Then there exists $\gamma>0$ depending only on $d$ and $V$ such that
$$
G(x,t)-p\cdot x\leq -\gamma t +C.
$$
\end{lem}

Proof: Since $|p|=1$, without loss of generality, we may assume 
that $p\cdot (1,0)>{1\over 2}$. Owing to the above Lemma \ref{lem:direction-reachability}, 
there exists $\beta>0$ such that, starting from any point $x\in Q=[0,\pi]^2$, player I can design a strategy so that it takes at most $2\beta $ time to move the game trajectory into $U_1+(-2\pi,0)$. Using periodicity, staring from any point $x\in Q$, in time $t$, player I can design a strategy to move the game trajectory into $U_1+(-2\pi k,0)$ for some $k\geq {t\over 2\beta}-1$ within time t.  Since $U_1+(-2\pi k,0)$ is $-V$ flow invariant,  when the game trajectory arrives at  $U_1+(-2\pi k,0)$, if needed, player I can just choose $\eta=0$ until time $t$. Accordingly, by the game theory formulation,  we have that 
$$
G(x,t)\leq \max_{y\in Q+k(-2\pi,0)} p\cdot y \leq -kp\cdot (2\pi,0)+C\leq -\gamma t +C.
$$
for $\gamma={\pi \over 2\beta}$. \qed

Combining with Lemma \ref{lem:ergodic-version},  we deduce that
\begin{cor}\label{cor:strictnegative} When $\lambda$ is small enough, 
$$
\max_{x\in \Rset^2}\lambda v_{\lambda}(x)\leq -{\gamma\over 2}.
$$
In particular, $v=v_{\lambda}$ is a viscosity supersolution of 
$$
\left(1-d\Div{\frac{p+Dv}{|p+Dv|}}\right)_{+}|p+Dv|+V(x)\cdot (p+Dv) 
={\gamma\over 2}>0, \quad \text{in $\Rset^2$}.
$$
Here $\gamma$ is from the previous Lemma \ref{lem:cauchy-limit}.
\end{cor}

Let $O_1=(0,0)$, $O_2=(1,0)$, $O_3=(1,1)$ and $O_4=(0,1)$. 
\begin{lem}\label{lem:boundness} 
Let $u_\lambda(x)=p\cdot x+v_{\lambda}(x)$. Then
$$
\max_{(x,y)\in [-\pi,\pi]^2}|u_{\lambda}(x)-u_{\lambda}(y)|\leq C
$$
for a constant $C$ depending only on $V$ and $d$.
\end{lem}

Proof: Owing to (\ref{v-max}), $u_\lambda$ is always a viscosity subsolution of
$$
\underline F(D^2u_\lambda,Du_\lambda)+V(x)\cdot Du_\lambda=1+\max_{\Rset^2}|V(x)| \quad \text{in $\Rset^2$}.
$$
Throughout the proof, $C_\mu$ represents various constants depending only on $d$, $V$ and a given parameter $\mu$.

\medskip

{\bf Step 1:} We first establish the difference bound within $Q_\mu$ for any given $\mu\in (0,1)$. Note that the level curve $\{x\in Q|\ H(x)=\mu\}$ is flow invariant.   When $\mu$ is close to 1, the traveling time around the level curve  is near $2\pi$. When $\mu$ is close to 0, the traveling time is bounded by $O(|\log\mu|)$.

Owing to Lemma \ref{lem:interior-reach}, Lemma \ref{lem:reachcontrol} and Lemma \ref{lem:flowcontrol}, we have that 

\medskip

1. For every point $x\in \partial Q_{{\mu}}$ and every point $y\in \overline Q_{{\mu}}$, $u_{\lambda}(x)\geq u_{\lambda}(y)-C_{\mu}$. Accordingly, 
$$
\min_{x\in \partial Q_{{\mu}}}u_{\lambda}(x)\geq \max_{x\in \overline Q_{{\mu}}}u_{\lambda}(x)-C_{\mu}.
$$

2. By Corollary \ref{cor:strictnegative} and the minimum value principle Lemma \ref{lem:mini},
$$
\min_{x\in \partial Q_{{\mu}}}u_{\lambda}(x)=\min_{x\in \overline Q_{{\mu}}}u_{\lambda}(x).
$$

Accordingly, 
\be\label{interior-control}
\max_{x\in \overline Q_{{\mu}}}u_{\lambda}(x)-\min_{x\in \overline Q_{{\mu}}}u_{\lambda}(x)=\max_{x,y\in \overline Q_{{\mu}}} |u_\lambda(x)-u_\lambda(y)|\leq C_{\mu}.
\ee

\medskip

{\bf Step 2:} Recall the four cells in $[-\pi,\pi]^2$: $U_1=(0, \pi)\times (0,\pi)$, $U_2=U_1-(\pi,0)$, $U_3=U_1-(0,\pi)$ and $U_4=U_1-(\pi,\pi)$. Let  $P_0=({\pi\over 2}, {\pi\over 2})$, $q_1=(0,0)$, $q_2=-(\pi,0)$, $q_3=-(0,\pi)$ and $q_4=-(\pi,\pi)$.  Then the above (\ref{interior-control}) also holds when $Q_\mu$ is replaced by $Q_{\mu}+q_i$ for $1\leq i\leq 4$. 

Let $\mu_0$ be the number from  Lemma \ref{lem:boundary-reach} that works for all four cells.  Combining with  corresponding versions of Corollary \ref{cor:cross-reach} in all four cells,   $Q_{2\mu_0}+q_i$ is reachable by $P_0+q_j$  for $1\leq i,j\leq 4$ within time $T_0$ depending only on $d$ and  $V$.  Since $Q_{2\mu_0}+q_i$ are flow invariant, owing to Lemma \ref{lem:reachcontrol} and (\ref{interior-control}) from Step 1, we deduce that for $1\leq i,j\leq 4$
\be\label{D-cell-com}
|u_\lambda(P_0+q_i)-u_{\lambda}(P_0+q_j)|\leq C. 
\ee
 By employing  (\ref{interior-control}) in all cells, we actually have that for any $\mu\in (0,1]$, 
\be\label{D-cell-com-1}
|u_\lambda(x)-u_{\lambda}(y)|\leq C_\mu  
\ee
for $x,y\in \cup_{1\leq i\leq 4}(Q_{\mu}+q_i)$. 
\medskip

{\bf Step 3:} It remains to take care of regions near  $\{H=0\}$ since $C_\mu\to +\infty$ as $\mu\to 0$.  Owing to   Corollary \ref{cor:cross-reach}, there exists $\tilde t_1$ such that every point $(1, d)$ for $d\in [0,1]$ is reachable from $P_0$ within time $\tilde t_1$.  

 For $d\in [0,1]$, let  $\xi_d:[0,\infty)\to \Rset^2$ be the $-V$ flow starting from $(1,d)$, i.e.,  
$\dot \xi_d(s)=-V(\xi_d(s))$ subject to $\xi_d(0)=(1,d)$.  
Write $\xi_0(\tilde t_1)=(\pi-2\nu_0, 0)$.  Fix  $0<\alpha_0\leq \min\{\mu_0, \nu_0\}$.  
 Clearly,  there exists $d_0\in (0, \mu_0)$ such that  for all $d\leq d_0$, 
$$
\xi_d([0,\infty))\cap (\{\pi-\alpha_0\}\times [0, 1])=\xi_d(s_d)
$$
for some $s_d\geq \tilde t_1$.  Let
$$
(\pi-\alpha_0, \beta_0)=\xi_{d_0}(s_{d_0}). 
$$
By choose $d_0$ small enough, we may assume that $\beta_0\leq \alpha_0$. 
\begin{center}
\includegraphics[scale=0.4]{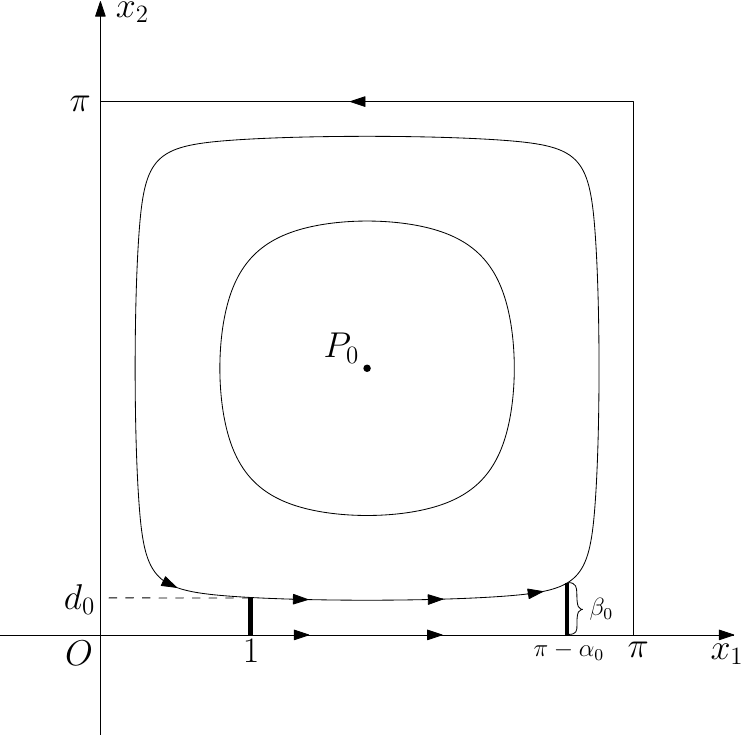}
\captionof{figure}{$d_0$, $\beta_0$ and $\alpha_0$. }
\end{center}
 Due to Lemma  \ref{lem:reachcontrol}, for all $d\leq d_0$,
$$
u_{\lambda}(P_0)\leq \max_{t\in [0, \tilde t_1]}u_\lambda(\xi_d(t))+C.
$$
Combining with Lemma \ref{lem:flowcontrol},  for all $s\in [0,1]$, 
$$
u_{\lambda}(P_0)\leq u_\lambda(\pi-\alpha_0, s\beta_0)+C.
$$
By looking at other cells and using (\ref{D-cell-com}),   we can find common $\beta_0\leq \alpha_0\in (0, \mu_0)$ such that
$$
u_{\lambda}(P_0)\leq \min\{u_{\lambda}(x)|\ x\in J_1\cup (-J_1)\cap J_2\cap (-J_2)\}+C.
$$
Here $J_1=\{\pi-\alpha_0\}\times [-\beta_0, \beta_0]$ and $J_2= [-\beta_0, \beta_0]\times \{\alpha_0\}$.  
\begin{center}
\includegraphics[scale=0.8]{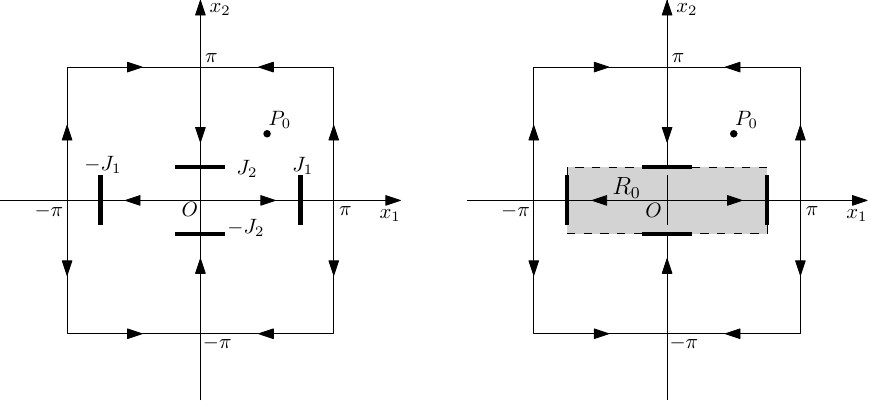}
\captionof{figure}{$J_i$ and $R_0$. }
\end{center}

Combining with Lemma \ref{lem:boundary-reach}, Lemma  \ref{lem:reachcontrol}, the flow-invariance of $Q_{\mu}$ and (\ref{D-cell-com-1}), we have that
$$
u_{\lambda}(P_0)\geq \max\{u_{\lambda}(x)|\ x\in J_1\cup (-J_1)\cap J_2\cap (-J_2)\}-C. 
$$
Now consider  the rectangle $R_0=[-\pi+\alpha_0, \pi-\alpha_0]\times [-\alpha_0, \alpha_0]$.  Thanks to (\ref{D-cell-com-1}) and the above argument,  we have that
$$
\max_{x\in \partial R_0}|u_{\lambda}(P_0)-u_\lambda(x)|\leq C.
$$
Again, combining with Lemma \ref{lem:boundary-reach}, Lemma  \ref{lem:reachcontrol}, the flow-invariance of $Q_{\mu}$ and (\ref{D-cell-com-1}), 
$$
u_\lambda(P_0)\geq  \max_{x\in R_0}u_\lambda(x)-C.
$$
Finally, applying the minimum value principle Lemma \ref{lem:mini} on $R_0$, we have that 
$$
\min_{x\in R_0}u_\lambda(x)=\min_{x\in \partial R_0}u_\lambda(x)\geq u_\lambda(P_0)-C.
$$
Accordingly, combining all the inequalities above, we derive that
$$
\max_{x\in R_0}|u_{\lambda}(P_0)-u_\lambda(x)|\leq C.
$$
By applying this to other cells and edges (see Fig. 12 below), we have that 
$$
\max_{x\in \Gamma_{\alpha_0}}|u_{\lambda}(P_0)-u_\lambda(x)|\leq C.
$$
Similar conclusion also holds in other cells. Combining with (\ref{D-cell-com-1}),  
 the lemma holds. \qed

\begin{center}
\includegraphics[scale=0.8]{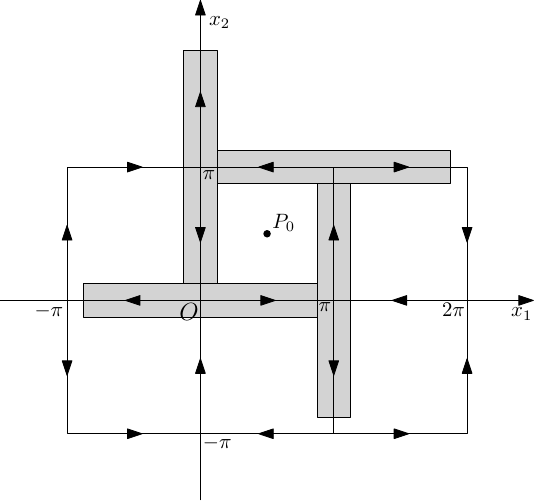}
\captionof{figure}{Corresponding versions of  $R_0$ on other edges. }
\end{center}

According to the above Lemma \ref{lem:boundness},  standard arguments lead to 
\begin{cor}\label{cor:unilimit} For any unit vector $p$, there exists $\overline H(p)\in \Rset$ such that
$$
\lim_{\lambda\to 0}\lambda v_{\lambda}(x)=-\overline H(p) \quad \text{uniformly in $\Rset^2$}.
$$
In particular, there exists a $\Zset^2$-periodic $\overline v\in USC(\Rset^2)$ that is a viscosity subsolution of 
$$
\underline F(D^2\overline v, p+D\overline v)+V(x)\cdot (p+D\overline v)=\overline H(p)
$$
and a $\Zset^2$-periodic $\underline v\in LSC(\Rset^2)$ that is a viscosity supersolution of 
$$
\overline F(D^2\underline v , p+D\underline v)+V(x)\cdot (p+D\underline v)=\overline H(p).
$$
Also,
$$
\sup_{x\in \Rset^2}\{|\overline {v}(x)|,\ |\underline{v}(x)|\}<C_0 
$$
for a constant $C_0$ depending only on $d$ and $V$. 
\end{cor}
Proof: The argument of the existence of the limit is similar to the proof of (3) in the following lemma. Let
$$
\overline v(x)=\limsup_{\lambda\to 0,y\to x}(v_{\lambda}(y)-v_{\lambda}(0))
$$
and
$$
\underline v(x)=\liminf_{\lambda\to 0,y\to x}(v_{\lambda}(y)-v_{\lambda}(0)).
$$
It is easy to see that our conclusion holds. 
\qed

\medskip

\begin{lem}\label{lem:conhbar} Below are several properties of $\overline H$. 

1. for $\lambda>0$, 
$$
\overline H(\lambda p)=\lambda \overline H(p);
$$

2. 
$$
\overline H(p)>0 \quad \text{for $p\not=0$};
$$

3.
$$
\overline H(p)\in C(\Rset^2).
$$
\end{lem}

Proof: (1) is obvious from the definition of $v_\lambda$. (2) follows from Corollary \ref{cor:strictnegative}. (3) follows from the standard stability of viscosity solutions. For the reader's convenience, we present details of the proof. Fixed a unit vector $p$. We argue by contradiction. If not, without loss of generality, we may assume that there exist a sequence of unit vectors $\{p_m\}_{m\geq 1}$ such that for $\lim_{m\to +\infty}p_m=p$ and for all $m\geq 1$
$$
\overline H(p_m)\geq \overline H(p)+r
$$
for some $r>0$. Owing to Corollary \ref{cor:unilimit}, for each $m\geq 1$, let $v_m\in LSC(\Rset^2)$ be a $\Zset^2$-periodic viscosity supersolution of 
$$
\overline F(D^2v_m, p_m+Dv_m)+V(x)\cdot (p_m+Dv_m)= \overline H(p)+r
$$
and
$$
\sup_{x\in \Rset^2} |v_m(x)|\leq C_0
$$
with the constant $C_0$ from Corollary \ref{cor:unilimit}.
Let 
$$
\underline v(x)=\liminf_{m\to +\infty, y\to x}v_{m}(y).
$$
Then $\underline v(x)\in LSC(\Rset^2)$ is a $\Zset^2$-periodic viscosity supersolution of 
$$
\overline F(D^2\underline{v}, p+D\underline {v})+V(x)\cdot (p+D\underline {v})= \overline H(p)+r.
$$
Also by Corollary \ref{cor:unilimit}, there is $\overline v\in USC(\Rset^2)$ that is a $\Zset^2$-periodic viscosity subsolution to
$$
\overline F(D^2{\overline v}, p+D {\overline v})+V(x)\cdot (p+D\overline {v})= \overline H(p).
$$
Write $\underline u=p\cdot x+\underline v$ and $\overline u=p\cdot x +\overline v$. Then $\underline u$ is a viscosity supersolution of 
$$
\overline F(D^2\underline u, D\underline u)+V(x)\cdot D\underline u=\overline H(p)+r
$$
and $\overline u$ is a viscosity subsolution of 
$$
\overline F(D^2 \overline u, D \overline u)+V(x)\cdot D\overline u=\overline H(p).
$$
Let
$$
w_{\delta}(x,y)=\overline u(x)-\underline u(y)-{|x-y|^4\over \delta}.
$$
Due to the periodicity of $\overline v$ and $\underline v$, we can find $(x_{\delta}, y_{\delta})\in \Rset^2\times \Rset^2$ such that
$$
w_{\delta}(x_{\delta}, y_{\delta})=\max_{x,y\in \Rset^2}w_{\delta}(x,y).
$$
It is easy to see that 
\be\label{zero-limit}
\lim_{\delta\to 0} w_{\delta}(x_{\delta}, y_{\delta})=\max_{x\in \Rset^2}(\overline v-\underline v)\quad \mathrm{and} \quad \lim_{\delta\to 0} {|x_{\delta}-y_{\delta}|^4\over \delta}=0.
\ee
Owing to Theorem 3.2 and Remark 3.8 in \cite{CL1992}, there are two $2\times 2$ symmetric matrices $X$ and $Y$ such that
$$
X\leq Y, \quad ||X||+||Y||\leq C{|x_{\delta}-y_{\delta}|^2\over \delta}
$$
and for $\bar p=4(x_{\delta}-y_{\delta}){|x_{\delta}-y_{\delta}|^2\over \delta}$
$$
\begin{array}{ll}
&\underline F(X,\bar p)+V(x_\delta)\cdot \bar p\leq \overline H(p)\\[3mm]
&\overline F(Y,\bar p)+V( y_\delta)\cdot \bar p\geq \overline H(p)+r.
\end{array}
$$
Also, 
$$
|(V(x_\delta)-V(y_\delta))\cdot \overline p|\leq {C|x_{\delta}-y_{\delta}|^4\over \delta}.
$$
Accordingly,
$$
\underline F(X,\bar p)-\overline F(Y,\bar p)\leq {C|x_{\delta}-y_{\delta}|^4\over \delta}-{r}.
$$
Case 1: If $x_\delta-y_{\delta}=0$, then $X=Y=0$. We have that
$$
0=0-0\leq -{r},
$$
which is impossible. 
\medskip

Case 2: If $x_\delta-y_{\delta}\not =0$. Then
$$
\underline F(X,\bar p)-\overline F(Y,\bar p)=F(X,\bar p)-F(Y,\bar p)\geq 0,
$$
we obtain
$$
0\leq {C|x_{\delta}-y_{\delta}|^4\over \delta}-{r}.
$$

\medskip

Owing to (\ref{zero-limit}), as $\delta\to 0$, we have
$$
0\leq -{r},
$$
a contradiction. \qed

\medskip

{\bf Proof of Theorem \ref{theo:main}} Let $\overline {v}$ and $\underline {v}$ be functions from Corollary \ref{cor:unilimit}. Then by comparison principle Theorem \ref{theo:comparison},
$$
p\cdot x+\epsilon \left( \overline v\left({x\over \epsilon}\right)-C_0\right)-\overline H(p)t\leq G_{\epsilon}(x,t)\leq p\cdot x+\epsilon \left(\underline v\left({x\over \epsilon}\right)+C_0\right)-\overline H(p)t
$$
and
$$
p\cdot x+\overline v(x)-C_0-\overline H(p)t\leq G(x,t)\leq p\cdot x+\underline v(x)+C_0-\overline H(p)t.
$$
for the constant $C_0$ from Corollary \ref{cor:unilimit}. The conclusion follows immediately. \qed

\begin{rmk}\label{rmk:property}  For application purposes, people are  interested in deriving explicit formulas of the effective burning velocity under the G-equation model  (see \cite{KAW} for instance).  Although a simple formula is mathematically not available,  it might be practically interesting to investigate more detailed properties of $\overline H(p)$ (e.g.,  its anisotropy  due to the presence of the fluid) in addition to its dependence on physical parameters. These kind of problems often require methods deeper than those standard PDE approaches in Lemma \ref{lem:conhbar}. For instance,    questions in this aspect have been studied in \cite{JTY} and \cite{TY} for the case  $v_{\vec n}=a(x)$ (e.g., phase transition in an inhomogeneous  medium without considering curvature effect) using tools from dynamical systems.  
A conclusion from results there is that in 2D, 
a polygon could be an effective front if and only if it is centrally symmetric with rational vertices and nonempty interior.  So different distributions of defects and heterogeneities can change the evolution significantly, which could  have an  important practical implication \cite{CB}.

\end{rmk}

\begin{rmk}\label{rmk:cm}  A natural question is whether our result can be 
proved by pure PDE or geometric approaches as in previous literature. 
For this aspect, Theorem 10.2 in \cite{CM} seems  relevant. 
Checking Assumption B'' there (if it holds) might need to evolve $Q_\mu$ by  
properly combining the motion law and the $V$ flow. It is not clear to us how 
to arrange the motions to reach a stationary supersolution (basically a large $V$-flow 
invariant set in our context) within finite time. Game theoretical method provides more 
flexibility to handle detailed local structures. For instance, the non-divergence 
inviscid example in section 11.2 of \cite{CM} relying on Theorem 10.2 for 2D can be proved 
for all dimensions via control formulation using full reachability or 
partially reachability $+$ minimum value principle as in this paper. 
See \cite{XY2014} for more applications of game theory in homogenization of 
non-coercive non-convex G-equations.
\end{rmk}

\section{Extension to general 2D incompressible flows}\label{sec:extension} 
In this section, we will briefly explain how to possibly modify our methods to 
cover more general 2D incompressible flows. Assume that $V=(-H_{x_2},H_{x_1})$ for a 
general periodic stream function $H$. For simplicity, we just discuss two representative scenarios
 and the corresponding adjustments of our methods. 

\medskip

{\bf Case 1:} Non-convex cells. In general, the level curves of $H$ might consist 
of convex and concave parts (relative to a fixed cell).
\begin{center}
\includegraphics[width=0.5\linewidth]{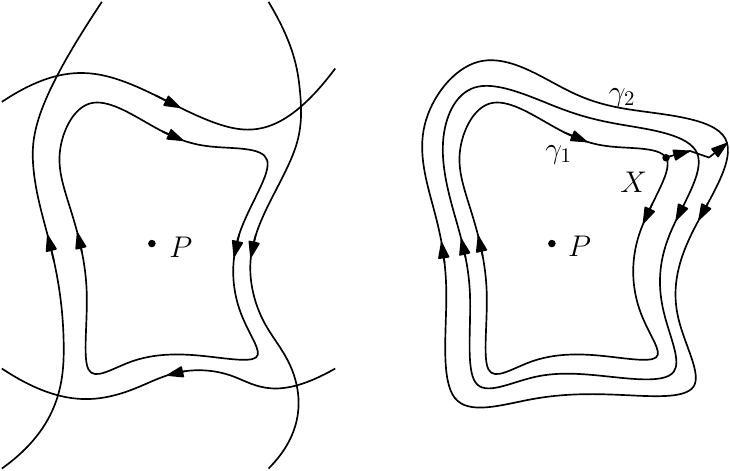}
\captionof{figure}{Move forward within a non-convex cell.}
\end{center}

$\bullet$ Movement within a cell (extension of Lemma \ref{lem:interior-reach}). 
Starting from a level curve $\gamma_1$, the game trajectory can move forward on the 
convex part (or near flat part) of $\gamma_1$ using a modified version of Lemma \ref{lem:perreach} 
to arrive at another level curve $\gamma_2$, then travel along $\gamma_2$ to its 
convex parts (or near flat parts) in order to further move forward. 

\medskip

$\bullet$ Movement to reach the boundary and adjacent cells 
(extension of Lemma \ref{lem:boundary-reach} and Corollary \ref{cor:cross-reach}). 
Starting from the center critical point $P$ in the cell $U_0$, player I aims to 
drive the game trajectory to travel to a point $W_2$ that is very close to a concave portion 
of the boundary. As in this paper, we need to control the amount of travel time. 
Since the game trajectory might not be able to cross the boundary on the concave portion, 
it can first move to a point $W_1$ near (but not very close to) the convex part of the boundary. 
Next use a modified version of Lemma \ref{lem:perreach} to reach $W_3$ in the adjacent cell $U_1$. 
Then follow the flow $-V$ to  a point $W_4$ that is close (but not too close) to the concave 
part of the boundary of $U_0$ and move back to $U_0$ based on a modified version of 
Lemma \ref{lem:perreach}, and finally travel to $W_2$ along the flow $-V$. 

Other lemmas and conclusions in this paper can be extended similarly. 
\begin{center}
\includegraphics[width=0.5\linewidth]{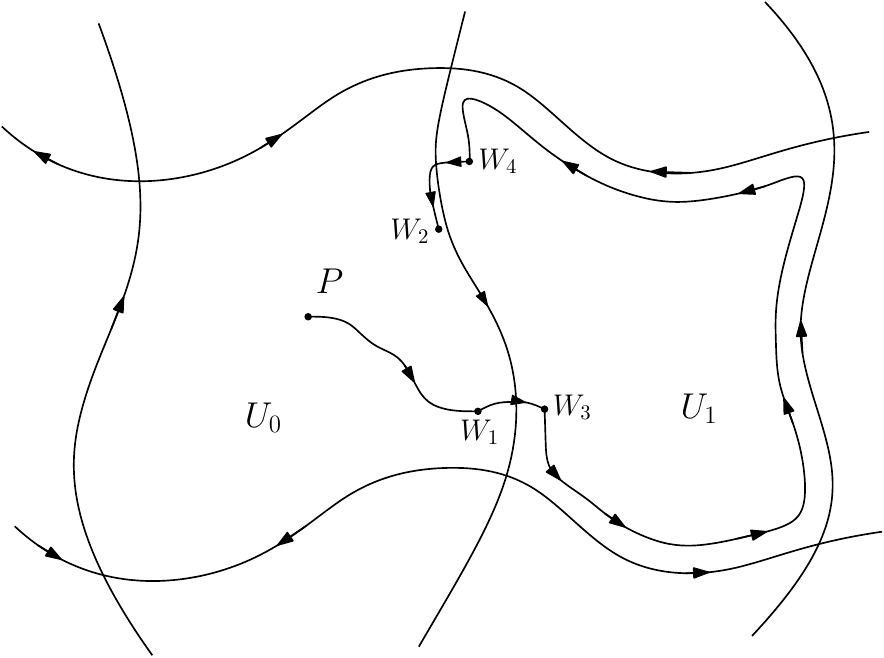}
\captionof{figure}{Move to points near the boundary and to an adjacent cell.}
\end{center}

\medskip

{\bf Case 2:} Cat's-eye types.  For simplicity,  we look at  the following  representative example
$$
H(x_1,x_2)=\sin x_1\sin x_2+\delta \cos x_1\cos x_2
$$
for $\delta\in (0,1)$. The picture consists of islands (e.g. shaded regions $I_1$, $I_2$ in figure below) and unbounded periodic orbits of $-V$ flow (e.g. regions $J_1$, $J_2$ in the figure below). 
\begin{center}
\includegraphics[width=0.5\linewidth]{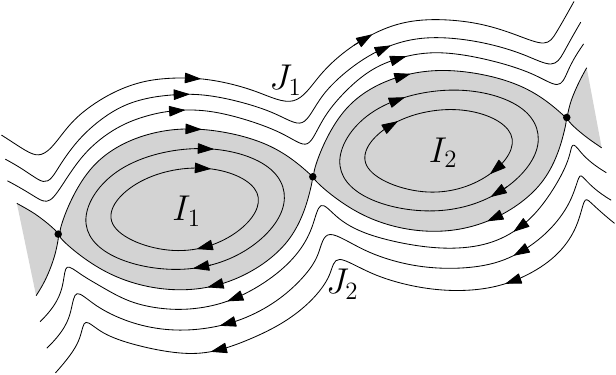}
\captionof{figure}{Cat's-eye type flows.}
\end{center}

Step 1: By similar arguments in this paper and possible extensions as in the above 
Case 1 together with periodicity, there exist constants $\bar c_1$, $\bar c_2$, $c_1$ 
and $c_2$ such that for $i=1,2$,  $\lim_{\lambda\to 0}\lambda u_{\lambda}=\bar c_i$ in 
the unbounded region $J_i$ and $\lim_{\lambda\to 0}\lambda u_{\lambda}=c_i$ within $I_i$.

\medskip

Step 2:  By partial reachability (from $I_1$, $I_2$ to $J_1$ and $J_2$ due to convexity 
of the boundary of $I_i$), we have that $c_1, c_2\leq \bar c_1, \bar c_2$.

\medskip

Step 3:  Note that for $(x_1,x_2)$ near $(0,0)$,
$$
H(x_1,x_2)=x_1x_2+\delta\left(1-{x_1^2\over 2}-{x_2^2\over 2}\right)+O(x_1^4+x_2^4).
$$
It is not hard to see that  the curvature of $H(x_1,x_2)=\delta$ (the boundary of islands)  tends 0 when approaching saddle points (inflection points).  Therefore,  starting from some regions of $J_1$ or $J_2$, player I can move the game trajectory into those islands through portions of the boundary near saddle points. Accordingly, we have that $c_1=c_2=\bar c_1=\bar c_2$. 

For other cat's-eye type flows, it could happen that one island $I$ is uniformly convex and, hence, 
  game trajectories might not be able to enter it from unbounded domains when the Markstein 
diffusivity $d$ is large. For this case, game trajectories can still enter islands adjacent 
to $I$ and the minimum value principle will lead to the same conclusion. One such example is 
$$
H(x_1,x_2)=x_1^2-x_2^2-2x_1^3+x_1^4  \quad \text{for $(x_1,x_2)$ near $(0,0)$}.
$$
Here  $H=0$ implies $x_2=x_1-x_1^2$ or $x_2=-x_1+x_1^2$ near $(0,0)$.   

 We expect that the effective burning velocity should exist for quite general 2D 
incompressible flows, at least if all critical points $H$ are non-degenerate where 
the flow structure is well-understood and essentially the combination of case 1 and 
case 2 (see \cite{A1991}).  The main challenge is how to find an efficient systematic 
proof without examining all possible scenarios. For non-smooth flows ($H\in  C^{1,1}$, or equivalently  
$V$ is only Lipschitz continuous) or flows with degenerate stagnation points, 
the analysis appears more complicated. We plan to investigate this issue in the future. 

\begin{rmk}\label{rmk:non-periodic} Our proof of Theorem 1.1. suggests a {\it general framework}  based on the one sided reachability of 
game trajectories and the minimum value principle of suitable stationary problems. 
The former is a separate dynamical system issue that has to be verified for a given flow field V, 
be it periodic, almost periodic or random. If it is true or almost surely true, the remaining PDE argument 
extends and the existence of  effective front speed might be established in interesting non-periodic settings. 
\end{rmk}

\section{Appendix} 
It is well known that the front propagation under the viscosity solution framework is consistent with the classical meaning when the smooth solutions exist. See section 6 in \cite{ES} for instance. The following conclusion is a 
special case in our context, which is needed to derive a reachability property. For the reader's convenience, 
we present its proof here in 2D, which is sufficient for our purpose. 

{Throughout this section,  we only assume that  $V\in W^{1,\infty}(\Rset^2)$ and 
\be\label{eq:vanishingline}
V(x)\cdot (1,0)=0 \quad \text{for $x\in \{0\}\times [0,1]$}.
\ee
} Let $S=(0,1)^2$ and 
$$
g_S(x)=
\begin{cases}
-{2\over \pi}\arctan(d(x,\partial S)) \quad \text{for $x\in S$}\\[3mm]
{2\over \pi}\arctan(d(x,\partial S)) \quad \text{otherwise}.
\end{cases}
$$
Here $d(x,\partial S)$ is the distance from $x$ to $\partial S$. 

\begin{lem}\label{lem:pass}
Suppose that $G\in C(\Rset^2\times [0,\infty))$ is the unique viscosity solution to equation (\ref{ge1}) subject to 
$$
G(x,0)=g_S(x).
$$
Then for a given $\delta\in (0,{1\over 2})$, there exists $t_{\delta}>0$ depending only on $d$, $V$ and $\delta$ such that 
$$
G((0,\theta), t)<0 \quad \text{for $(\theta, t)\in [\delta, 1-\delta]\times (0, t_\delta]$}. 
$$
\end{lem}

Proof: Intuitively, this conclusion is obvious since the speed along the normal direction $\vec{n}=(-1,0)$ at the point $(0, \theta)$ is 
$$
v_{\vec{n}}=1-d\kappa+V(x)\cdot \vec{n}=1.
$$
To make this rigorous, we need to build smooth supersolutions and employ comparison principle. It suffices to prove this at a fixed $\theta\in [\delta,1-\delta]$. 

\medskip

{\bf Step 1: Choose an ellipse.} Consider the ellipse 
\begin{equation*}
E_{\theta}(t): \frac{(x_1-a_0 - \nu)^2}{a^2(t)} + \frac{(x_2-\theta)^2}{b^2(t)} = 1.
\end{equation*}
Here $\nu\in (0,a_0)$ is added for technical convenience and will be sent to zero later. 

\medskip

\begin{center}
\includegraphics[width=0.4\linewidth]{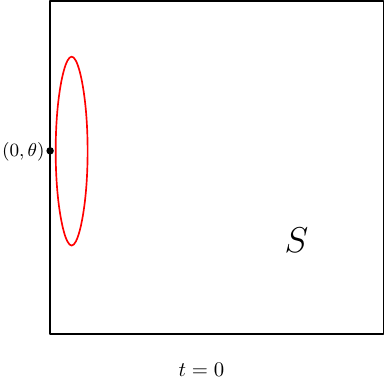}\qquad 
\includegraphics[width=0.45\linewidth]{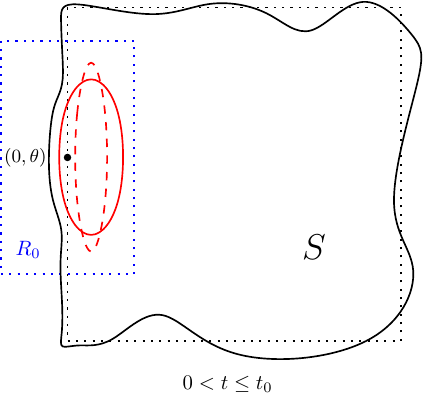}
\captionof{figure}{Propagation of the ellipse.}
\end{center}

Let $b_0={\delta\over 2}$. Then choose $a_0\in (0, {b_0\over 4})$ small enough such that

$$ 
|V(x)\cdot (1,0)|<{1\over 8} \quad \text{if $x\in [-4a_0, 4a_0]\times [0,1]$}
$$
and
$$
\frac{d64a_0}{b_{0}^2}+4M_0\sqrt{3}\frac{a_0}{b_{0}}<{1\over 8}.
$$
Here $M_0=\max_{x\in \Rset^2}|V(x)|$. 

Then we define $a(t) = a_0 + \frac{1}{2}t$ and $b(t) = b_{0} - Lt$ for $L>0$ satisfying  
$$
{1\over 2}-{3La_0\over 4b_0}<1-{db_0\over {a_0}^{2}}-M_0.
$$
Hereafter we require $0 \leq t \leq t_\delta$ for $t_\delta = \min \left\lbrace 2a_0, \frac{b_{0}}{2L}\right\rbrace$, which implies that
$$
(a(t), b(t))\in [a_0,\ 2a_0]\times \left[ \frac{b_0}{2},\ b_0\right]
$$ 
and
$$
E_\theta(t)\subset R_0=[-4a_0, 4a_0]\times [0,1]
$$
for $t\in [0,t_\delta]$. Also, $4a_0<b_0$ leads to 
$$
a(t)< b(t) \quad \text{for all $t\in [0,t_\delta]$}.
$$

\medskip

For convenience, we drop the dependence of $a$ and $b$ on $t$. Hereafter $t\in [0,t_\delta]$ unless 
specified otherwise. Let $\vec{n}$ be the outward unit normal vector along $E_{\theta}(t)$ that has the following parameterization: for $\phi\in [0, 2\pi]$
$$
\begin{cases}
x_1=a_0+\nu+a\cos \phi\\
x_2=\theta+b\sin \phi.
\end{cases}
$$
Then
\begin{equation*}
{V} \cdot \vec{n} = {V} \cdot \left( \frac{b\cos\phi}{\sqrt{a^2\sin^2\phi + b^2\cos^2\phi}}, \frac{a\sin\phi}{\sqrt{a^2\sin^2\phi + b^2\cos^2\phi}}\right). 
\end{equation*}

If $|\sin \phi|<{\sqrt{3}\over 2}$,
\be\label{v-bound}
\begin{array}{ll}
|{V} \cdot \vec{n}| &\leq |V(x)\cdot (1,0)| + M_0 \cdot \frac{a\sin\phi}{\sqrt{a^2\sin^2\phi + b^2\cos^2\phi}}\\[3mm]
&\leq |V(x)\cdot (1,0)| + M_0 \cdot \frac{a}{b} \cdot |\tan\phi|\\[3mm]
& \leq |V(x)\cdot (1,0)| + M_0 \cdot \frac{2a_0}{b_0/2} \cdot \sqrt{3} \\[3mm]
&= |V(x)\cdot (1,0)|+ 4\sqrt{3}M_0 \cdot \frac{a_0}{b_0}< {1\over 8}+{1\over 8}={1\over 4}.
\end{array}
\ee

{\bf Step 2: Evolution of an elliptic boundary.} Let us recall some basic facts. Given a $C^1$ function $f(x,t)$ and the family of level curves 
$$
C(t)=\{x\in \Rset^2\ | f(x,t)=0\}, 
$$
if $D_xf\not =0$, the propagation speed of $C(t)$ along the outward normal direction $\vec{n}={D_xf\over |D_xf|}$ is given by
$$
v_{\vec n}=-{f_t\over |D_xf|},
$$
which can be easily derived through the chain rule. Moreover, the corresponding mean curvature along $C(t)$ is 
$$
\kappa=\mathrm{div}_x(\vec{n}).
$$

Now let us verify that for $t\in [0,t_\delta]$, the propagation of $E_{\theta}(t)$ obeys the 
following inequality:
\be\label{curve-motion}
v_{\vec n}< 1 - d\kappa + {V} \cdot \vec{n},
\ee
which will be used to construct a supersolution. Fix any $a, b > 0$ and an ellipse
\begin{equation*}
\frac{x_1^2}{a^2} + \frac{x_2^2}{b^2} = 1,
\end{equation*}
direct compuations show that its curvature at point $P(a\cos\phi, b\sin\phi)$, $0 \leq \phi < 2\pi$,  is 
\begin{equation*}
\kappa = \frac{ab}{\left(a^2\sin^2\phi + b^2\cos^2\phi\right)^{\frac{3}{2}}}
\end{equation*}
and the normal velocity at $P(a\cos\phi, b\sin\phi)$ is
\begin{eqnarray*}
v_{\vec{n}} = \frac{a^{\prime}b\cos^2\phi + ab^{\prime}\sin^2\phi}{\sqrt{a^2\sin^2\phi + b^2\cos^2\phi}}.
\end{eqnarray*}

\medskip

\noindent {\bf Case 1.} If $|\sin \phi|\leq {\sqrt{3}\over 2}$, then
\begin{equation*}
d\kappa = \frac{dab}{\left( a^2\sin^2\phi + b^2\cos^2\phi\right)^{\frac{3}{2}}} \leq \frac{dab}{\left( b^2\cos^2\phi \right)^{\frac{3}{2}}} \leq \frac{dab}{(b/2)^3}\leq \frac{d64a_0}{b_{0}^2} <{1\over 8}.
\end{equation*}

Note that
$$
v_{\vec n}= \frac{\frac{b}{2}\cos^2\phi - La \sin^2\phi}{\sqrt{a^2 \sin^2\phi + b^2\cos^2\phi}} \leq {{b\over 2}\cos^2\phi\over \sqrt{b^2\cos^2\phi}}\leq \frac{1}{2}.
$$ 
Since for $t\in [0,t_\delta]$, 
$$
E_{\theta}(t)\subset [-4a_0, 4a_0]\times [0,1],
$$
(\ref{v-bound}) implies that $v_{\vec{n}} \leq {1\over 2}< 1-{1\over 8}-{1\over 4}< 1 -d \kappa + {V} \cdot \vec{n}$.

\medskip 

\noindent {\bf Case 2.} If $|\sin \phi|\geq {\sqrt{3}\over 2}$, then
\begin{eqnarray*}
v_{\vec n}= \frac{\frac{b}{2}\cos^2\phi - La \sin^2\phi}{\sqrt{a^2 \sin^2\phi + b^2\cos^2\phi}} &\leq& {{b\over 2}\cos^2\phi\over \sqrt{b^2\cos^2\phi}}-{La\sin^2\phi\over b}\\
&\leq& {1\over 2}-{3La\over 4b}\leq {1\over 2}-{3La_0\over 4b_0}.
\end{eqnarray*}

Meanwhile,
$$
d\kappa = \frac{dab}{\left(a^2\sin^2\phi + b^2\cos^2\phi\right)^{\frac{3}{2}}}\leq {db\over a^2}\leq {db_0\over {a_0}^{2}}.
$$

Therefore, due to the choice of $L$, we have $v_{\vec n}< 1 - d\kappa + {V} \cdot \vec{n}$.

Combining case 1 and case 2, we see that (\ref{curve-motion}) holds, i.e.,  the evolution of $E_{\theta}(t)$ satisfies 
\begin{eqnarray*}
v_{\vec n} < 1 -d \kappa + {V} \cdot \vec{n}, &\text{for}& 0 \leq t \leq t_\delta. 
\end{eqnarray*}

\medskip

{\bf Step 3: Comparison.} Let
$$
\Psi(x,t)=\frac{(x_1-a_0 - \nu)^2}{a^2(t)} + \frac{(x_2-\theta)^2}{b^2(t)} -1.
$$

Owing to (\ref{curve-motion}), we may choose $\mu_0\in (0, {1\over 2})$ such that (\ref{curve-motion}) also holds along the curve $\{x\in \Rset^2|\ \Psi(x,t)=\mu\}$ for all $\mu\in [-\mu_0,\mu_0]$ and $t\in [0,t_\delta]$. Equivalently, for 
$$
D_{\theta}=\left\{(x,t)\in \Rset^2\times [0,t_\delta]\ |-\mu_0\leq \Psi(x,t) \leq \mu_0\right\},
$$
\be
\Psi_t + \left(1-d\, \, \mathrm{div}\left({D\Psi\over |D\Psi|}\right)\right)|D\Psi|+V(x)\cdot D\Psi\geq 0 \quad \text{on $D_{\theta}$}.
\ee
Let $\{h_k\}_{k\geq 1}\in C^{\infty}(\Rset)$ be a sequence of functions such that
$$
0<h_{k}^{'}\leq 1 \quad \text{in $I_k=\left(-\mu_0+{1\over k},\ \mu_0-{1\over k}\right)$},\quad h_{k}^{'}=0 \quad \text{in $\Rset\backslash I_k$} 
$$
and $\lim_{k\to +\infty}h_k(s)=h(s)$ uniformly in $\Rset$, where
$$
h(s)=
\begin{cases}
\mu_0 \quad \text{for $s\geq \mu_0$}\\
s\quad \text{for $s\in[-\mu_0, \mu_0]$}\\
-\mu_0 \quad \text{for $s\leq -\mu_0$}.
\end{cases}
$$
Apparently, $\Psi_k(x,t)=h_k(\Psi(x,t))$ satisfies 
\be
{\partial \Psi_{k}\over \partial t} + \left(1-d\, \, \mathrm{div}\left({D\Psi_k\over |D\Psi_k|}\right)\right)|D\Psi_k|+V(x)\cdot D\Psi_k\geq 0 \quad \text{on $\Rset^2\times (0,t_\delta)$}.
\ee
By stability, we have that $G_1(x,t)=h(\Psi(x,t))\in W^{1,\infty}(\Rset^2\times [0,\infty))$ is a viscosity supersolution of 
\be
{\partial G_{1}\over \partial t }+ \left(1-d\, \, \mathrm{div}\left({DG_1\over |DG_1|}\right)\right)|DG_1|+V(x)\cdot DG_1\geq 0 \quad \text{on $\Rset^2\times (0,t_\delta)$}.
\ee

Since $(a)_+\geq a$, $G_1=G_1(x,t)$ is also a viscosity supersolution of equation (\ref{ge1}) on $\Rset^2\times (0,t_\delta)$. 

Because for fixed $\nu>0$
$$
\{G_1(x,0)\leq 0\}=\{\Psi(x,0)\leq 0\}\subset S,
$$
we can choose a function $\xi\in C^{\infty}(\Rset)$ such that $\dot \xi>0$, $\xi(0)=0$, $\sup_{s\in \Rset}|\xi(s)|<\infty$ and 
$$
g(x)\leq \xi(G_1(x,0)).
$$
Since $\xi(G_1(x,0))$ is also a viscosity supersolution of equation (\ref{ge1}) on $\Rset^2\times (0,t_\delta)$, thanks to Theorem \ref{theo:comparison}, we have that
$$
G(x,t)\leq \xi(G_1(x,t)) \quad \text{for $(x,t)\in \Rset^2\times [0,t_\delta]$}.
$$
In particular, this implies that
$$
\{x\in \Rset^2|\ G_1(x,t)<0\}=\{x\in \Rset^2|\ \xi(G_1(x,t))<0\}\subset \{x\in \Rset^2|\ G(x,t)<0\}.
$$
Note that $\{x\in \Rset^2|\ G_1(x,t)<0\}=\{x\in \Rset^2|\ \Psi(x,t)<0\}$. Sending $\nu\to 0$, we have that for $t\in [0,t_\delta]$, 
$$
\left\{x=(x_1,x_2)\in \Rset^2|\ \frac{(x_1-a_0)^2}{a^2(t)} + \frac{(x_2-\theta)^2}{b^2(t)} <1\right\}\subset \{x\in \Rset^2|\ G(x,t)<0\}.
$$
Then
$$
G((0,\theta) ,t)<0 \quad \text{for $t\in (0, t_\delta]$},
$$
which finishes the proof. Note that $t_\delta$ only depends on $\delta$ and $V$. \qed

Let $\Omega\subset \Rset^2$ be an open convex set. Denote by $G_{\Omega}(x,t)$ 
the unique viscosity solution to equation (\ref{ge1}) subject to $G_\Omega(x,0)=g_{\Omega}(x)$ where 
$$
g_{\Omega}(x)=
\begin{cases}
-{2\over \pi}\arctan(d(x, \partial \Omega)) \quad \text{for $x\in \Omega$}\\[3mm]
{2\over \pi}\arctan(d(x,\partial \Omega) )\quad \text{otherwise}.
\end{cases}
$$
Given two sets $E_1$ and $E_2$, their Hausdorff distance 
$$
d_H(E_1, E_2)=\max\{ \max_{x\in E_1}d(x, E_2), \ \max_{x\in E_2}d(x, E_1)\}.
$$
Also, for $\alpha>0$ and $\delta\in (0, {1\over 2})$, we write
$$
W_{\alpha, \delta}=[-\alpha, \alpha]\times [\delta, 1-\delta].
$$

\begin{lem}\label{lem:stability}
Let $S=(0,1)^2$. For given $\delta\in (0,{1\over 2})$ and $n\in \Nset$, there exists $\sigma_{\delta, n}>0$ such that if 
$$
d_H(S, \Omega)\leq \sigma_{\delta,n} \quad \mathrm{and} \quad \alpha\leq \sigma_{\delta,n},
$$
then
$$
G_\Omega(x,t)<0 \quad \text{for $(x,t)\in W_{\alpha,\delta}\times \left[{t_{\delta}\over n},\ t_{\delta}\right]$}.
$$
Here $t_\delta$ is from the previous Lemma \ref{lem:pass}. 
\end{lem}
Proof: We argue by contradiction. If not, then there exist a sequence of convex open sets $\{\Omega_m\}_{m\geq 1}$ such that
$$
d_H(S, \Omega_m)\leq {1\over m}
$$
and for some $(x_m,t_m)\in W_{{1\over m},\delta}\times \left[{t_{\delta}\over n},\ t_{\delta}\right]$
$$
G_{\Omega_m}(x_m,t_m)\geq 0.
$$
Since 
$$
\lim_{m\to +\infty}g_{\Omega_m}(x)=g_S(x) \quad \text{uniformly on $\Rset^2$},
$$
due to Remark \ref{rmk:boundary}, the uniqueness of viscosity solutions, we have that 
$$
\lim_{m\to +\infty}G_{\Omega_m}(x,t)=G(x,t) \quad \text{locally uniformly on $\Rset^2\times [0,\infty)$}.
$$
Here $G$ is from Lemma \ref{lem:pass}. The proof is similar to that of (\ref{limit}). Also, up to a subsequence if necessary, we may assume that $\lim_{m\to +\infty}(x_m,t_m)=((0,\theta), \bar t)$ for $(\theta, \bar t)\in [\delta, 1-\delta]\times \left[{t_{\delta}\over n},\ t_{\delta}\right]$. Then we have that
$$
G((0,\theta), \bar t)=\lim_{m\to +\infty}G_{\Omega_m}(x_m,t_m)\geq 0.
$$
This is a contradiction. \qed

As an immediate corollary, we have the following reachability.

\begin{lem}\label{lem:perreach} Consider the game in section 2.2. Under the assumption of Lemma \ref{lem:stability}, every point on $W_{\alpha,\delta}$ can reach $\Omega$ within time $t_{\delta}\over n$. Also, it is easy to see that if we replace $S$ by an arbitrary rectangle, all previous results are still true in the corresponding forms.  See the figure below.\end{lem}
\begin{center}
\includegraphics[width=0.5\linewidth]{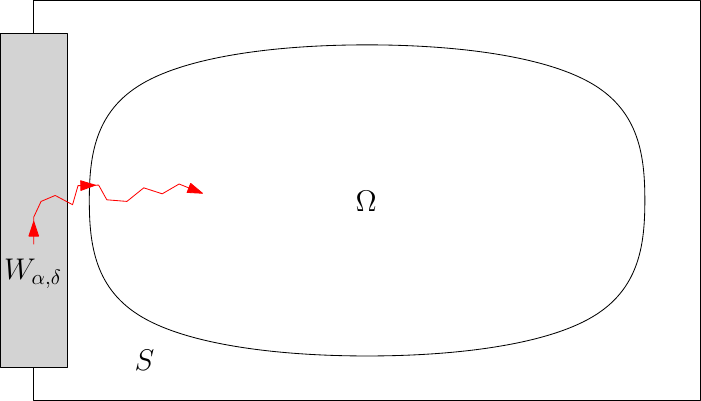}
\captionof{figure}{Reach $\Omega$ from $W_{\alpha,\delta}$. }
\end{center}

\begin{rmk}\label{rmk:randomwalk} Due to the hidden stochastic nature of the game trajectory, 
it is not clear to us how to use pure game dynamics to prove the above reachability conclusion. 
An interesting analog in random walk (or Brownian motion) is to use strong maximum principle of the 
Laplace equation to show that a particle has a positive probability 
to exit from any small window of the boundary. 
\end{rmk}

\medskip

{\bf Data availability statement.} This paper does not generate any data.

\medskip

{\bf Acknowledgements.} The authors would like to thank 
Inwon Kim and Hung V. Tran for  helpful comments in 
improving the presentation of the paper. 
We are grateful to Robert V. Kohn for fruitful discussions 
of game representations and for providing references. 
Last but not the least, we thank 
the anonymous referees for all the constructive comments. 

\bibliographystyle{plain}

\end{document}